\documentclass[a4paper,12pt]{article}

\usepackage{amsthm}
\usepackage{amsmath,amssymb,latexsym,amsfonts,mathrsfs}

\theoremstyle{plain}
\newtheorem{Thm}{Theorem}[section]
\newtheorem{Lem}[Thm]{Lemma}
\newtheorem{Prop}[Thm]{Proposition}
\theoremstyle{definition}

\newtheorem{Rem}[Thm]{Remark}
\newcommand{\Proof}[2][Proof]{\begin{proof}[{#1}] #2 \end{proof}}

\renewcommand{\d}{{\rm d}} %differential
\newcommand{\e}{{\rm e}} 
\newcommand{\eps}{\ensuremath{\varepsilon}}

\newcommand{\n}{\nonumber}

\numberwithin{equation}{section}

\makeatletter
\renewcommand\section{\@startsection {section}{1}{\z@}%
                                   {-3.5ex \@plus -1ex \@minus -.2ex}%
                                   {2.3ex \@plus.2ex}%
                                   {\normalfont\large\bf}}
\makeatother

\makeatletter
\renewcommand\subsection{\@startsection {subsection}{1}{\z@}%
                                   {-3.5ex \@plus -1ex \@minus -.2ex}%
                                   {2.3ex \@plus.2ex}%
                                   {\normalfont\normalsize\bf}}
\makeatother

\makeatletter
\@addtoreset{footnote}{page}
\makeatother

\newcommand{\absol}[1]{\left| #1 \right|} %absolute value
\newcommand{\rbra}[1]{\!\left( #1 \right)} %round brackets or parentheses
\newcommand{\cbra}[1]{\!\left\{ #1 \right\}} %curly brackets or braces
\newcommand{\sbra}[1]{\!\left[ #1 \right]} %square brackets or brackets
\newcommand{\pmat}[1]{\begin{pmatrix} #1 \end{pmatrix}} %()

\newcommand{\cB}{\ensuremath{\mathcal{B}}}
\newcommand{\cD}{\ensuremath{\mathcal{D}}}
\newcommand{\cG}{\ensuremath{\mathcal{G}}}
\newcommand{\cL}{\ensuremath{\mathcal{L}}}

\begin{document}

\begin{center}
{\Large \bf 
Extensions of diffusion processes on intervals and Feller's boundary conditions 
}
\end{center}
\begin{center}
Kouji \textsc{Yano}\footnote{
Graduate School of Science, Kyoto University, Kyoto, Japan.}\footnote{
Research partially supported by KAKENHI (20740060)}
%Firstname \textsc{Surname}\footnote{
%Affiliation
%}
\end{center}
\begin{center}
{\small \today}
\end{center}
\bigskip

%\begin{center}
%{\small Dedicated to      on the occasion of his        birthday}
%\end{center}
%\bigskip

\begin{abstract}
For a minimal diffusion process on $ (a,b) $, 
any possible extension of it to a standard process on $ [a,b] $ 
is characterized by 
the characteristic measures of excursions away from the boundary points $ a $ and $ b $. 
The generator of the extension is proved to be characterized by Feller's boundary condition. 
\end{abstract}


\section{Introduction}

The generator of a minimal diffusion process on an interval $ (a,b) $ 
can be characterized by the second order differential operator 
$ \cL = \cD_m \cD_s $ in Feller's canonical form, 
where $ m $ and $ s $ are strictly-increasing continuous functions on $ (a,b) $. 
Feller (\cite{MR0047886} and \cite{MR0092046}) determined 
all possibilities of boundary conditions 
to make $ \cL $ a generator of some Feller semigroup, 
which we call {\em Feller's boundary conditions}. 
If the boundary points $ a $ and $ b $ are both accessible, 
Feller's boundary condition is of the form 
\begin{align}
\Phi_a(f) = \Phi_b(f) = 0 
\label{eq: Feller bc}
\end{align}
where 
\begin{align}
\Phi_a(f) 
= p_1 f(a) - p_2 \cD_sf(a) + p_3 \cL f(a) - p_4[f-f(a)] 
\label{eq: Phia}
\end{align}
and 
\begin{align}
\Phi_b(f) 
= q_1 f(b) + q_2 \cD_sf(b) + q_3 \cL f(b) - q_4[f-f(b)] 
\label{eq: Phib}
\end{align}
for some non-negative constants $ p_i,q_i $, $ i=1,2,3 $ 
and some non-negative measures $ p_4 $ on $ (a,b] $ and $ q_4 $ on $ [a,b) $. 
Here, for a measure $ \mu $ and an integrable function $ f $, 
we write $ \mu[f] = \int f \d \mu $ . 

The purpose of this paper is to determine, for a given minimal diffusion process on $ (a,b) $, 
all possibilities of its extensions to standard processes on $ [a,b] $. 
For such an extension, 
we give the characteristic measure of excursions away from the boundary $ \{ a,b \} $, 
construct the process by piecing together excursions 
and prove that its $ C_b $-generator is characterized by Feller's boundary condition. 
The results of this paper complement the ealier results 
obtained by It\^o--McKean \cite{MR0154338}, It\^o \cite{ItoKyoto} and Rogers \cite{MR701528}. 
For an invariance principle for such processes, see Yano \cite{MR2543582}. 
Note that Fukushima \cite{FukuBdry} 
recently determined, by the help of the Dirichlet form theory, 
all possibilities of the extensions 
to symmetric diffusion processes 
for any minimal one-dimensinoal diffusion process possibly with killing inside $ (a,b) $. 
Our objects exclude killing inside $ (a,b) $, 
but are out of scope of \cite{FukuBdry} 
since they have jumps from the boundary so that they are not necessarily symmetric.

For this purpose, we shall appeal to the excursion theory; 
the sample path which starts from $ a $ 
up to the killing time or the first hitting time of $ b $, whichever comes earlier, 
will be constructed from the stagnancy rate $ p_3 $ 
and the characteristic measure of excursions away from $ a $ given as 
\begin{align}
n^{(b)}_a = p_1 \delta_{\Delta} + p_2 n^{(b)}_{a,{\rm refl}} 
+ \int_{(a,b]} p_4(\d x) P^{\rm stop}_x , 
\label{}
\end{align}
where $ \Delta $ stands for the excursion which stays at the cemetery for all time, 
$ n^{(b)}_{a,{\rm refl}} $ for the characteristic measure of excursions away from $ a $ 
of the diffusion process reflected at $ a $ and stopped at $ b $, 
and $ P^{\rm stop}_x $ for the law of the $ \cL $-diffusion process 
started at $ x $ and stopped at $ a $ or $ b $. 
The possible behaviors at $ a $ of a generic sample path are as follows: 
\begin{enumerate}
\item 
it is killed (i.e., jumps to the cemetery) according to the rate $ p_1 $; 
\item 
it enters the interior $ (a,b) $ continuously according to the rate $ p_2 $; 
\item 
it stagnates at $ a $ according to the rate $ p_3 $; 
\item 
it jumps into an interior point $ x $ of $ (a,b] $ 
picked accordingly to the rate $ p_4(\d x) $. 
\end{enumerate}
We may carry out the same construction of 
the sample path which starts from $ b $ as that from $ a $, 
where $ p_i $'s are replaced by $ q_i $'s. 
We shall prove that 
the $ C_b $-generator of the process constructed in such a way as above 
is characterized by Feller's boundary conditions \eqref{eq: Feller bc}. 

This paper is organized as follows. 
In Section \ref{sec: main}, we introduce several notations from the theory of excursions 
and state our main theorems. 
We give a very general formula about the resolvents 
and give theorems about $ C_b $-generators. 
In Section \ref{sec: resolv}, we prove the resolvent formula. 
In Section \ref{sec: minimal}, we study several properties 
of the resolvent operator of the minimal diffusion process. 
Section \ref{sec: determ} 
is devoted to the proofs of the theorems of $ C_b $-generators.

\section{Notations and main theorems} \label{sec: main} 

\subsection{Excursions and resolvents} \label{sec: eandr}

Let $ \Delta $ denote an isolated point attached to $ [-\infty ,\infty ] $ 
and we call $ \Delta $ the {\em cemetery}. 
For any function $ f $ defined on an subinterval of $ [-\infty ,\infty ] $, 
we always extend $ f $ so that $ f(\Delta) = 0 $. 
Let $ W $ denote the space consisting of all c\`adl\`ag functions 
$ w:[0,\infty ) \to [-\infty ,\infty ] \cup \{ \Delta \} $ such that, 
for some $ 0 \le T_{\Delta} \le \infty $, 
$ w(t) \in [-\infty ,\infty ] $ for all $ 0 \le t < T_{\Delta} $ 
and $ w(t) = \Delta $ for all $ t \ge T_{\Delta} $. 
We equip $ W $ with the $ \sigma $-field generated by the cylinder sets. 
We denote the coordinate process of $ W $ by $ (X(t))_{t \ge 0} = (X_t)_{t \ge 0} $: 
\begin{align}
X(t)(w) = X_t(w) = w(t) 
\quad \text{for} \ t \ge 0 . 
\label{}
\end{align}
We denote the first hitting time of $ x \in [-\infty ,\infty ] \cup \{ \Delta \} $ by 
\begin{align}
T_x = T_x(X) = \inf \{ t>0 : X_t = x \} 
\label{}
\end{align}
and the first exit time of $ x \in [-\infty ,\infty ] $ by 
\begin{align}
\tau_x = \tau_x(X) = \inf \{ t>0 : X_t \neq x \} 
\label{}
\end{align}
Here we adopt the usual convention that $ \inf \emptyset = \infty $.

Let $ a \in [-\infty ,\infty ] $. 
We write $ E(a) $ for the set of all $ e \in W $ 
such that 
$ e(t) $ is constant for all $ t \ge T_a \wedge T_{\Delta} $, 
where $ s \wedge t = \min \{ s,t \} $. 
Each element $ e $ of $ E(a) $ is called an {\em excursion away from $ a $}. 
We utilize the theory of excursions away from $ a $; 
see, e.g., \cite[Chap. III]{Blu} and \cite[Chap. IV]{Ber}. 
Let $ I $ be a subinterval of $ [-\infty ,\infty ] $ 
and let $ \{ (X_t)_{t \ge 0},(P_x)_{x \in I} \} $ 
be a standard process on $ I $, i.e., 
a strong Markov process such that 
$ X_t \in I $ for all $ t<T_{\Delta} $ 
a.s. with respect to $ P_x $ for $ x \in I $, 
and such that 
its sample paths are right-continuous 
and quasi-left-continuous up to $ T_{\Delta} $; 
see, e.g., \cite{BG}. 
Suppose that $ a \in I $. 
We write $ P^{(a)}_x $ for the law under $ P_x $ 
of the stopped process $ (X^{(a)}_t)_{t \ge 0} $ defined as 
\begin{align}
X^{(a)}_t = 
\begin{cases}
X_t & \text{if $ t < T_a $}, \\
a & \text{if $ t \ge T_a $}. 
\end{cases}
\label{}
\end{align}

For the process $ \{ (X_t)_{t \ge 0},P_a \} $, 
we would like to obtain a positive additive functional $ (L_a(t))_{t \ge 0} $ 
which increases only when the process stays at $ a $, 
and a non-negative constant $ \varsigma_a $ called the {\em stagnancy rate}. 
The following three cases are only possible: 
\begin{enumerate}
\item 
Suppose that $ a $ is {\em regular-for-itself}, i.e., $ P_a(T_a=0)=1 $, 
and also that $ a $ is {\em instantaneous}, i.e., $ P_a(\tau_a=0)=1 $. 
Let $ L_a(t) $ be a choice of the local time at $ a $ up to time $ t $. 
Then there exists a constant $ \varsigma_a \ge 0 $ such that 
\begin{align}
\int_0^t 1_{\{ X_s=a \}} \d s = \varsigma_a L_a(t) . 
\label{}
\end{align}
\item 
Suppose that $ a $ is regular-for-itself but not instantaneous. 
Then $ a $ is {\em holding}, i.e., $ P_x(\tau_a>0)=1 $. 
Let $ \varsigma_a>0 $ be a fixed constant and set 
\begin{align}
L_a(t) = \frac{1}{\varsigma_a} \int_0^t 1_{\{ X_s=a \}} \d s . 
\label{}
\end{align}
\item 
Suppose that $ a $ is irregular-for-itself. Then we have $ P_a(T_a>0)=1 $. 
Then the set $ \{ t>0:X_t=a \} $ is locally finite 
so that can be enumerated as $ \{ \tau_1 < \tau_2 < \cdots \} $. 
Extending the probability space, 
we let $ \{ e_n \}_{n=0}^{\infty } $ be a sequence of independent standard exponential variables 
which is independent of $ (X_t)_{t \ge 0} $. 
We now define 
\begin{align}
L_a(t) = e_0 + \sum_{n=1}^{\infty } e_n 1_{\{ \tau_n \le t \}} 
\label{}
\end{align}
and define $ \varsigma_a=0 $. 
\end{enumerate}

Let $ \eta_a(l) = \inf \{ t>0: L_a(t)>l \} $. 
Then, in any one of the above three cases (i)-(iii), 
the process $ (\eta_a(l))_{l \ge 0} $ is a subordinator 
which explodes at $ L_a(\infty ) $. 
We define a point process $ (p_a(l))_{l \in D_a} $ taking values in $ E(a) $ as 
\begin{align}
p_a(l)(t) = 
\begin{cases}
X(t+\eta_a(l-)) & \text{if $ 0 \le t < \eta_a(l) - \eta_a(l-) $}, \\
X(\eta_a(l)) & \text{if $ t \ge \eta_a(l) - \eta_a(l-) $} 
\end{cases}
\label{}
\end{align}
for all $ l $ belonging to the domain $ D_a = \{ l : \eta_a(l-)<\eta_a(l) \} $. 
We see that $ (p_a(l))_{l \in D_a} $ is a stationary Poisson point process 
stopped at $ L_a(\infty ) $, 
so that it is called the {\em excursion point process}. 
Its characteristic measure will be denoted by $ n_a $ 
and called the {\em characteristic measure of excursions away from $ a $}. 
The process $ (\eta_a(l))_{l \ge 0} $ can be recovered from the excursion point process as 
\begin{align}
\eta_a(l) = \varsigma_a l + \sum_{s \le l} T_a(p_a(s)) . 
\label{eq: eta-a(l)}
\end{align}
By \eqref{eq: eta-a(l)}, we have $ P_a[\e^{-r \eta_a(l)}] = \e^{-l \psi_a(r)} $, where 
\begin{align}
\psi_a(r) = \varsigma_a r + n_a \sbra{ 1-\e^{-r T_a} } . 
\label{}
\end{align}
Since $ P_a(\eta_a(l)<\infty ) > 0 $, we obtain 
\begin{align}
n_a \sbra{ 1-\e^{- T_a} } < \infty , 
\label{eq: cond1}
\end{align}
and, in particular, $ n_a $ is $ \sigma $-finite. 
If $ a $ is regular-for-itself, 
the process $ \eta_a(l) $ is strictly increasing until it explodes, 
so that we have 
\begin{align}
\varsigma_a > 0 
\quad \text{or} \quad 
[n_a(T_a<t) = \infty \ \text{for all $ t>0 $}]. 
\label{eq: cond2}
\end{align}
We know that the strong Markov property holds in the following sense: 
for any stopping time $ T $, any constant time $ t $ 
and any non-negative Borel functions $ f $ and $ g $, one has 
\begin{align}
n_a \sbra{ f(X_T) g(X_{T+t}) } = n_a \sbra{ f(X_T) P^{(a)}_{X_T}[g(X_t)] } . 
\label{}
\end{align}
If we write 
\begin{align}
k_a = n_a(\{ \Delta \}) 
, \quad 
\nu_a(\cdot) = n_a(\cdot;X_0=a) 
, \quad 
j_a(\cdot) = n_a(X_0 \in \cdot) , 
\label{}
\end{align}
then $ n_a $ may be represented as 
\begin{align}
n_a = k_a \delta_{\Delta} + \nu_a 
+ \int_{I \setminus \{ a \}} j_a(\d x) P^{(a)}_x . 
\label{eq: repre of na}
\end{align}
The quantities $ k_a $, $ \nu_a $ and $ j_a $, respectively, 
are called the {\em killing rate}, 
the {\em characteristic measure of excursions of continuous entrance}, 
and the {\em jumping-in measure}, respectively. 

\begin{Rem} \label{rem: irregular}
Suppose that $ 0 < n_a(X_0=a) < \infty $. 
Then the first index $ \lambda $ such that $ p_a(\lambda)(0)=a $ 
is positive and finite a.s., so that 
$ \eta_a(\lambda) $ is a finite stopping time such that $ X(\eta_a(\lambda)) = a $. 
Hence it follows that $ a $ may not be regular-for-itself 
and that $ n_a = \nu_a = P^{(a)}_a $. 
\end{Rem}

Let $ \cB_b(I) $ or simply $ \cB_b $ 
denote the set of all bounded Borel functions on $ I $. 
Let us study the resolvent operator: 
\begin{align}
R_rg(x) = P_x \sbra{ \int_0^{\infty } \e^{-rt} g(X_t) \d t } 
\quad \text{for $ g \in \cB_b $, $ r>0 $ and $ x \in I $}. 
\label{eq: ro}
\end{align}
Let $ b \in I \setminus \{ a \} $ be fixed and write 
\begin{align}
T_{a,b} = T_a \wedge T_b = \inf \{ t>0 : X_t \in \{ a,b \} \} . 
\label{}
\end{align}
Let $ \varsigma^{(b)}_a $ and $ n^{(b)}_a $ be the stagnancy rate 
and the characteristic measure of excursions away from $ a $ 
for the process $ \{ (X^{(b)}_t)_{t \ge 0},P_a \} $ 
starting from $ a $ and stopped at $ b $. 
For $ r>0 $, we define 
\begin{align}
\psi_a^{(b)}(r) = \varsigma_a^{(b)} r + n_a^{(b)} \sbra{ 1-\e^{-rT_a} } 
\label{}
\end{align}
and, for $ r>0 $ and $ g \in \cB_b $, we define 
\begin{align}
N^{(b)}_{a,r}(g) 
= \varsigma^{(b)}_a g(a) + n^{(b)}_a \sbra{ \int_0^{T_{a,b}} \e^{-rt} g(X_t) \d t } . 
\label{}
\end{align}
The following resolvent formula will play a key role: 

\begin{Thm} \label{thm: main-1}
One has 
\begin{align}
\psi_a^{(b)}(r) R_rg(a) 
= N^{(b)}_{a,r}(g) + n^{(b)}_a \sbra{ \e^{-r T_b} ; T_b < \infty } R_rg(b) . 
\label{}
\end{align}
\end{Thm}

Theorem \ref{thm: main-1} will be proved in Section \ref{sec: resolv}. 
We remark that Theorem \ref{thm: main-1} contains 
Theorem 5.3 of It\^o \cite{ItoKyoto} and Theorem 1 of Rogers \cite{MR701528} 
as special cases where $ b $ is not accessible. 

Before closing this subsection, we introduce another notation. 
Recall that the resolvent operator is defined in \eqref{eq: ro}. 
It is easy to see that the resolvent equation 
\begin{align}
R_r-R_q + (r-q) R_r R_q = 0 
\label{}
\end{align}
holds. 
Let $ C_b(I) $ or simply $ C_b $ denote the set of all bounded continuous functions on $ I $. 
If $ g \in C_b $, we see, by the dominated convergence theorem, 
that $ r R_r g(x) $ converges to $ g(x) $ as $ r \to \infty $ for all $ x \in I $. 
Hence we see that the range $ R_r(C_b) $ does not depend on $ r>0 $ 
and that $ R_r:C_b \to R_r(C_b) $ is injective. 
Now we define the operator $ \cG $ as 
\begin{align}
\begin{cases}
D(\cG) = R_r(C_b) , \\
\cG f = rf-R_r^{-1}f \quad \text{for $ f \in D(\cG) $} , 
\end{cases}
\label{}
\end{align}
which turns out to be independent of the choice of $ r>0 $. 
We call $ \cG $ the {\em $ C_b(I) $-generator} 
or simply the {\em $ C_b $-generator}.

\subsection{Feller's boundary classification}

Let us prepare several notations of Feller's characteristics of one-dimensional diffusions. 

Let $ (a,b) $ be a subinterval of $ [-\infty ,\infty ] $. 
Let $ m $ and $ s $ be strictly-increasing continuous functions 
on $ (a,b) $ with values in $ [-\infty ,\infty ] $. 
We extend such functions so that 
$ m(a)=m(a+) $, $ m(b) = m(b-) $, $ s(a)=s(a+) $ and $ s(b)=s(b-) $.  
For $ a \le x < y \le b $, we write $ m(x,y) = m(y)-m(x) $ and $ s(x,y) = s(y)-s(x) $. 

We adopt {\em Feller's classification of boundaries} 
taken from \cite{MR0047886} and \cite{MR0063607} 
as follows\footnote{The term ``enterable" in the second line of \eqref{eq: Feller1} 
has never been used in the literatures. 
We prefer, however, to adopt this unfamiliar terminology 
rather than to use the term ``entrance", which was used in It\^o--McKean's textbook \cite{IM}, 
for the sake of consistency with \eqref{eq: Feller2}. }. 
Let $ c \in (a,b) $. 
We use the following terminology: 
\begin{align}
& \begin{cases}
\text{$ a $ is called {\em accessible} if} \ \int_a^c \d s(x) \int_x^c \d m(y) < \infty , 
\\
\text{$ a $ is called {\em enterable} if} \ \int_a^c \d m(x) \int_x^c \d s(y) < \infty , 
\end{cases}
\label{eq: Feller1}
\end{align}
where we write $ \int_{c_1}^{c_2} $ for $ \int_{(c_1,c_2]} $ if $ c_1 \le c_2 $ 
and for $ - \int_{(c_2,c_1]} $ if $ c_1 > c_2 $. 
These classifications do not depend on the particular choice of $ c \in (a,b) $. 
We also use the following terminology: 
\begin{align}
& \begin{cases}
\text{$ a $ is called {\em regular} if it is both accessible and enterable} 
, \\
\text{$ a $ is called {\em exit} if it is accessible but not enterable} 
, \\
\text{$ a $ is called {\em entrance} if it is enterable but not accessible} 
, \\
\text{$ a $ is called {\em natural} if it is neither accessible nor enterable} . 
\end{cases}
\label{eq: Feller2}
\end{align}
The same classification is made for the boundary $ b $ 
by switching the roles between $ a $ and $ b $.

\subsection{Minimal one-dimensional diffusion processes} \label{sec: mp}

We call $ \{ (X_t)_{t \ge 0},(P^0_x)_{x \in (a,b)} \} $ a {\em minimal diffusion process} 
if it is a standard process on $ (a,b) $ such that the following conditions hold 
$ P_x $-a.s. for all $ x \in (a,b) $: 
\begin{enumerate}
\item 
$ t \mapsto X_t $ is continuous in $ t<T_{\Delta} $; 
\item 
$ X_t \in (a,b) $ for all $ t<T_{\Delta} $; 
\item 
$ X_t \to a $ or $ X_t \to b $ as $ t \to T_{\Delta}- $. 
\end{enumerate}
In addition, we assume the following condition: 
\begin{align}
P_x(T_y<\infty )>0 
\quad \text{for all $ x,y \in (a,b) $}. 
\label{}
\end{align}

We denote its resolvent by 
\begin{align}
R^0_rg(x) = P^0_x \sbra{ \int_0^{\infty } \e^{-rt} g(X_t) \d t } 
\quad \text{for $ g \in \cB_b((a,b)) $, $ r>0 $ and $ x \in (a,b) $}. 
\label{eq: resolv of mdp}
\end{align}
It is well-known that there exist some strictly-increasing continuous functions $ m $ and $ s $ 
such that 
\begin{align}
P^0_x(T_{a'}>T_{b'}) = \frac{s(a',x)}{s(a',b')} 
\quad \text{for all $ a < a' < x < b' < b $} 
\label{}
\end{align}
and 
\begin{align}
P^0_x[T_{a'} \wedge T_{b'}] = \int_{(a',b')} K_{a',b'}(x,y) \d m(y) 
\quad \text{for all $ a < a' < x < b' < b $} , 
\label{}
\end{align}
where 
\begin{align}
K_{a',b'}(x,y) = K_{a',b'}(y,x) = \frac{s(a',x) s(y,b')}{s(a',b')} 
\quad \text{for all $ a' \le x \le y \le b' $} . 
\label{}
\end{align}
The functions $ m $ and $ s $ will be called 
the {\em speed measure} and the {\em canonical scale}, respectively. 
We write $ D(\cD_s) $ for the set of measurable functions $ f $ on $ [a,b] $ such that 
$ f(x)-f(y) = \int_x^y g(z) \d s(z) $ for all $ x,y \in (a,b) $ 
for some locally $ s $-integrable function $ g $, 
and in this case, we write $ \cD_s f = g $. 
We write $ D(\cL) $ for the set of measurable functions $ f $ on $ [a,b] $ such that 
$ f \in D(\cD_s) $ and 
$ \cD_s f(x)-\cD_s f(y) = \int_x^y g(z) \d m(z) $ for all $ x,y \in (a,b) $ 
for some locally $ m $-integrable function $ g $, 
and in this case, we write $ \cL f = \cD_m \cD_s f = g $. 
For $ f \in D(\cL) $, we say that $ \cL f \in C_b([a,b)) $ 
if $ \cL f(a+) $ exists, and in this case we write $ \cL f(a) = \cL f(a+) $. 
We will use the following property for the resolvent $ (R^0_r)_{r>0} $. 

\begin{Prop} \label{prop: L0R}
For all $ g \in \cB_b((a,b)) $ and $ r>0 $, it follows that 
\begin{align}
\text{$ R^0_rg \in D(\cL) $ and $ \cL R^0_rg = r R^0_rg - g $} 
\label{}
\end{align}
and that $ f=R^0_rg $ satisfies the following: 
\begin{enumerate}
\item 
$ f(a+)=0 $, if $ a $ is accessible; 
\item 
$ f(a+) = u_r(a) \int_a^{b-} g(x) v_r(x) \d m(x) $ 
and 
$ \cD_s f(a+) = 0 $, if $ a $ is entrance; 
\item 
$ f(a+) = g(a+)/r $, if $ a $ is natural and the limit $ g(a+) $ exists. 
\end{enumerate}
\end{Prop}

For the proof of Proposition \ref{prop: L0R}, see, e.g., \cite[Theorem 62.1]{ItoEssential}. 
We note that, when $ a $ is entrance, 
$ r \int_a^{b-} v_r(x) \d m(x) = \cD_s v_r(b) - \cD_s v_r(a) < \infty $. 

For $ x=a $ and $ b $, let $ \gamma(x)=1 $ if $ x $ is accessible and $ \gamma(x)=0 $ otherwise. 
The $ C_b((a,b)) $-generator $ \cG^0 $ of the minimal diffusion process 
can be characterized as follows.

\begin{Thm} \label{thm: g of mp}
The $ C_b((a,b)) $-generator $ \cG^0 $ of the resolvent $ (R^0_r)_{r>0} $ is given as 
\begin{align}
D(\cG^0) = \{ f \in D(\cL) : f,\cL f \in C_b((a,b)) , \ \gamma(a)f(a+)=\gamma(b)f(b-)=0 \} 
\label{}
\end{align}
and 
\begin{align}
\cG^0 f(x) = \cL f(x) 
\quad \text{for} \ f \in D(\cG^0) \ \text{and} \ x \in (a,b). 
\label{}
\end{align}
\end{Thm}

Theorem \ref{thm: g of mp} is well-known, so that we omit the proof.

\subsection{Extensions of the minimal process} \label{sec: emp}

Let $ \{ (X_t)_{t \ge 0},(P^0_x)_{x \in (a,b)} \} $ be a minimal diffusion process 
and let $ m $ and $ s $ be the corresponding speed measure and the canonical scale, respectively. 
Let $ I = [a,b) $ or $ [a,b] $ 
and let $ \{ (X_t)_{t \ge 0},(P_x)_{x \in I} \} $ 
be an extension to a standard process of the minimal process 
$ \{ (X_t)_{t \ge 0},(P^0_x)_{x \in (a,b)} \} $, i.e., 
the killed process $ (X^0_t)_{t \ge 0} $ defined by 
\begin{align}
X^0_t = 
\begin{cases}
X_t & \text{if $ 0 \le t < T_a \wedge T_b $}, \\
\Delta & \text{if $ t \ge T_a \wedge T_b $} 
\end{cases}
\label{}
\end{align}
considered under $ P_x $ has the law $ P^0_x $ for all $ x \in (a,b) $. 
We have essentially the following four cases: 
\begin{description}
\item[{\boldmath $ 1^{\circ}). $}] 
$ I = [a,b) $, $ a $ is accessible and $ b $ is not; 
\item[{\boldmath $ 2^{\circ}). $}] 
$ I = [a,b] $, $ a $ is accessible, and $ b $ is not; 
\item[{\boldmath $ 3^{\circ}). $}] 
$ I = [a,b] $, and neither $ a $ nor $ b $ is accessible; 
\item[{\boldmath $ 4^{\circ}). $}] 
$ I = [a,b] $, and both $ a $ and $ b $ are accessible. 
\end{description}
In each case, 
we would like to characterize the $ C_b(I) $-generator of the extension. 

As we will discuss the case of {\boldmath $ 4^{\circ}) $} after a while, 
let us assume that $ a $ is accessible. 
For the stopped process $ \{ (X^{(b)}_t)_{t \ge 0},P_a \} $, 
we let $ \varsigma^{(b)}_a $ and $ n^{(b)}_a $ be 
the stagnancy rate and the characteristic measure of excursions away from $ a $. 
Let $ \varsigma^{(a)}_b $ and $ n^{(a)}_b $ be 
their counterparts for the stopped process $ \{ (X^{(a)}_t)_{t \ge 0},P_b \} $. 

If $ a $ is regular, then 
there exists, uniquely in the sense of law, a conservative diffusion process on $ [a,b] $ 
which is stopped at $ b $ 
and whose stagnancy rate at $ a $ is zero. 
This process will be called the {\em reflected process} at $ a $. 
We denote by $ n^{(b)}_{a,{\rm refl}} $ 
the characteristic measure of excursions away from $ a $ which is normalized so that 
\begin{align}
n^{(b)}_{a,{\rm refl}} (T_x<\infty ) = \frac{1}{s(a,x)} 
\quad \text{for} \ x \in (a,b) . 
\label{eq: nrefl}
\end{align}
We do not define $ n^{(b)}_{a,{\rm refl}} $ if $ a $ is exit. 

For the stopped process $ \{ (X^{(b)}_t)_{t \ge 0},P_a \} $, 
we have 
\begin{align}
n^{(b)}_a = p_1 \delta_{\Delta} + p_2 n^{(b)}_{a,{\rm refl}} 
+ \int_{(a,b]} p_4(\d x) P^{\rm stop}_x 
\label{eq: n(b)a}
\end{align}
for some non-negative constants $ p_1 $ and $ p_2 $ and 
some non-negative measure $ p_4 $; 
in fact, we have the representation \eqref{eq: repre of na} 
so that $ \nu_a $ should be proportional to $ n^{(b)}_{a,{\rm refl}} $ if $ a $ is regular. 
We let $ p_3 = \varsigma^{(b)}_a $. 
If $ a $ is exit, $ \nu_a $ should be zero, 
so that we let $ p_2=0 $ and discard the term $ p_2 n^{(b)}_{a,{\rm refl}} $. 
By the conditions \eqref{eq: cond1} and \eqref{eq: cond2}, 
the coefficients must satisfy the conditions: 
\begin{align}
\int_{(a,b]} p_4(\d x) P^{\rm stop}_x \sbra{ 1-\e^{-T_a} } 
< \infty 
\label{eq: pcond1}
\end{align}
and 
\begin{align}
p_2 + p_3 > 0 
\quad \text{or} \quad 
[p_4((a,a+\eps)) = \infty \ \text{for all $ \eps>0 $}]. 
\label{eq: pcond2}
\end{align}

\begin{Rem}
The condition \eqref{eq: pcond1} is equivalent to the following condition: 
\begin{enumerate}
\item 
$ \int_a^{a+\eps} p_4(\d x) s(a,x) < \infty $ for some $ \eps>0 $ 
if $ a $ is regular; 
\item 
$ \int_a^{a+\eps} p_4(\d x) \int_a^x m(y,c) \d s(y) < \infty $ for some $ \eps>0 $ 
if $ a $ is exit. 
\end{enumerate}
\end{Rem}

\begin{Rem}
The behavior at $ a $ of the process $ \{ (X_t)_{t \ge 0},(P_x)_{x \in [a,b]} \} $ is as follows: 
\begin{enumerate}
\item $ a $ is regular-for-itself and instantaneous 
if either one of the following holds: 
\\ \quad 
(i-1) $ a $ is regular and either $ p_2>0 $ or $ p_4((a,a+\eps)) = \infty $ for all $ \eps>0 $; 
\\ \quad 
(i-2) $ a $ is exit and $ p_4(a,a+\eps) = \infty $ for all $ \eps>0 $; 
\item $ a $ is regular-for-itself but not instantaneous 
if either one of the following holds: 
\\ \quad 
(ii-1) $ a $ is regular, $ p_2=0 $ and $ p_4((a,a+\eps)) < \infty $ for some $ \eps>0 $; 
\\ \quad 
(ii-2) $ a $ is exit and $ p_4((a,a+\eps)) < \infty $ for some $ \eps>0 $. 
\end{enumerate}
\end{Rem}

We state our main theorems which determine the $ C_b $-generator of all possible extensions. 
The proofs will be given in Section \ref{sec: determ}. 

{\boldmath $ 1^{\circ}). $} 
Suppose that $ I=[a,b) $, $ a $ is accessible and $ b $ is not. 
The following theorem generalizes Theorem 5.3 of It\^o \cite{ItoKyoto}. 

\begin{Thm} \label{thm: main--}
The $ C_b([a,b)) $-generator $ \cG $ is such that its domain is given as 
\begin{align}
D(\cG) = \cbra{ f \in D(\cL) : 
f, \cL f \in C_b([a,b)) \ \text{and} \ 
\Phi_a(f) = 0 } , 
\label{eq: DG domain--}
\end{align}
where $ \Phi_a $ has been defined in \eqref{eq: Phia}, 
and it satisfies 
\begin{align}
\cG f(x) = \cL f(x) 
\quad \text{for} \ f \in D(\cG) \ \text{and} \ x \in [a,b) . 
\label{}
\end{align}
\end{Thm}

\begin{Rem}
If, in addition, $ b $ is entrance, then one has, for any $ f \in D(\cG) $, 
\begin{align}
\cD_s f(b-) = 0 . 
\label{eq: dsfb}
\end{align}
This remark holds true also in any other case below. 
\end{Rem}

{\boldmath $ 2^{\circ}). $} 
Suppose that $ I=[a,b] $, $ a $ is accessible and that $ b $ is not. 
If $ b $ is entrance and irregular-for-itself, then we have 
\begin{align}
\text{$ \varsigma^{(a)}_b=0 $ and $ n^{(a)}_b = P^{\rm stop}_b $}. 
\label{}
\end{align}
Otherwise, we have 
\begin{align}
\text{$ q_3:=\varsigma^{(a)}_b>0 $ and 
$ n^{(a)}_b = q_1 \delta_{\{ \Delta \}} + \int_{[a,b)} q_4(\d x) P^{\rm stop}_x $} 
\label{}
\end{align}
for some non-negative constant $ q_1 $ and some non-negative finite measure $ q_4 $ on $ [a,b) $, 
and we put $ q_2=0 $ for convenience.

\begin{Thm} \label{thm: main-nonacc}
The domain of the $ C_b([a,b]) $-generator $ \cG $ is given as follows: 
if $ b $ is entrance and irregular-for-itself, then 
\begin{align}
D(\cG) = 
\{ f \in D(\cL) : f,\cL f \in C_b([a,b]), \ \Phi_a(f)=p_4(\{ b \}) f(b) \} ; 
\label{eq: DG domain1}
\end{align}
if $ b $ is natural or [entrance and regular-for-itself], then 
\begin{align}
D(\cG) = 
\{ f \in D(\cL) : f,\cL f \in C_b([a,b]), \ \Phi_a(f)=p_4(\{ b \}) f(b) , \ \Phi_b(f)=0 \} , 
\label{eq: DG domain2}
\end{align}
where $ \Phi_b $ has been defined in \eqref{eq: Phib}. 
In both cases, one has, for any $ f \in D(\cG) $, 
\begin{align}
\cG f(x) = \cL f(x) 
\quad \text{for} \ x \in [a,b] . 
\label{eq: G1}
\end{align}
\end{Thm}

{\boldmath $ 3^{\circ}). $} 
Suppose that $ I=[a,b] $ and neither $ a $ nor $ b $ is accessible. 
We only state a special case; 
the necessary modifications in the other cases are obvious. 
Let us assume that $ a $ is natural 
and that $ b $ is entrance and irregular-for-itself. 
For the stopped process $ \{ (X^{(b)}_t)_{t \ge 0},P_a \} $, 
we have 
\begin{align}
\text{$ p_3:=\varsigma^{(b)}_a>0 $ and 
$ n^{(b)}_a = p_1 \delta_{\{ \Delta \}} + \int_{(a,b]} p_4(\d x) P^{\rm stop}_x $} 
\label{}
\end{align}
for some non-negative constant $ p_1 $ and some non-negative finite measure $ p_4 $ on $ (a,b] $, 
and we put $ p_2=0 $ for convenience. 
For the stopped process $ \{ (X^{(a)}_t)_{t \ge 0},P_b \} $, 
we have 
\begin{align}
\text{$ \varsigma^{(a)}_b=0 $ and $ n^{(a)}_b = P^{\rm stop}_b $}. 
\label{}
\end{align}

\begin{Thm} \label{thm: main-nonacc2}
The $ C_b([a,b]) $-generator $ \cG $ 
is such that its domain is given as 
\begin{align}
D(\cG) =& 
\{ f \in D(\cL) : f,\cL f \in C_b([a,b]), \ \Phi_a(f)=p_4(\{ b \}) f(b) \} . 
\label{eq: DG domain2+}
\end{align}
For any $ f \in D(\cG) $, one has 
\begin{align}
\cG f(x) = \cL f(x) 
\quad \text{for} \ x \in [a,b] . 
\label{eq: G2}
\end{align}
\end{Thm}

{\boldmath $ 4^{\circ}). $} 
Suppose that $ I=[a,b] $ and both $ a $ and $ b $ are accessible. 
For the stopped process $ \{ (X^{(a)}_t)_{t \ge 0},P_b \} $, 
we have 
\begin{align}
n^{(a)}_b = q_1 \delta_{\Delta} + q_2 n^{(a)}_{b,{\rm refl}} 
+ \int_{[a,b)} q_4(\d x) P^{\rm stop}_x 
\label{eq: n(a)b}
\end{align}
for some non-negative constants $ q_1 $ and $ q_2 $ and 
some non-negative measure $ q_4 $ on $ [a,b) $ such that, 
\begin{align}
\int_{[a,b)} q_4(\d x) P^{\rm stop}_x \sbra{ 1-\e^{-T_a} } 
< \infty 
\label{eq: qcond1}
\end{align}
and, for $ q_3 := \varsigma^{(a)}_b $, 
\begin{align}
q_2 + q_3 > 0 
\quad \text{or} \quad 
[q_4((b-\eps,b)) = \infty \ \text{for all $ \eps>0 $}]. 
\label{eq: qcond2}
\end{align}

\begin{Thm} \label{thm: main}
The $ C_b([a,b]) $-generator $ \cG $ is such that its domain is given as 
\begin{align}
D(\cG) = \cbra{ f \in D(\cL) : 
f, \cL f \in C_b([a,b]) \ \text{and} \ 
\Phi_a(f) = \Phi_b(f) = 0 } , 
\label{eq: DG domain}
\end{align}
and it satisfies 
\begin{align}
\cG f(x) = \cL f(x) 
\quad \text{for} \ f \in D(\cG) \ \text{and} \ x \in [a,b] . 
\label{}
\end{align}
\end{Thm}

\section{The resolvent formula} \label{sec: resolv} 

Let $ \{ (X_t)_{t \ge 0},(P_x)_{x \in [a,b]} \} $ 
be a standard process and follow the notations in Subsection \ref{sec: eandr}. 
We utilize the following lemma: 

\begin{Lem} \label{lem: stagnancy}
For any $ r>0 $, one has 
\begin{align}
\int_0^{\eta_a(l)} \e^{-rt} 1_{\{ X(t)=a \}} \d t 
= \varsigma_a \int_0^l \e^{-r\eta_a(s)} \d s 
\quad \text{for all $ l \ge 0 $, $ P_a $-a.s.} 
\label{}
\end{align}
\end{Lem}

\Proof{
Take a sample point from a set of full probability. Note that 
\begin{align}
\int_0^{\eta_a(l)} 1_{\{ X(t)=a \}} \d t 
= \eta_a(l) - \sum_{s \le l} \{ \eta_a(s)-\eta_a(s-) \} 
= \varsigma_a l . 
\label{}
\end{align}
Let $ \eps>0 $. Then we have 
\begin{align}
\int_{\eta_a(l)}^{\eta_a(l+\eps)} \e^{-rt} 1_{\{ X(t)=a \}} \d t 
\le \e^{- r \eta_a(l)} \int_{\eta_a(l)}^{\eta_a(l+\eps)} 1_{\{ X(t)=a \}} \d t 
= \e^{- r \eta_a(l)} \varsigma_a \eps 
\label{eq: lem stag1}
\end{align}
and 
\begin{align}
\int_{\eta_a(l)}^{\eta_a(l+\eps)} \e^{-rt} 1_{\{ X(t)=a \}} \d t 
\ge \e^{-r\eta_a(l+\eps)} \int_{\eta_a(l)}^{\eta_a(l+\eps)} 1_{\{ X(t)=a \}} \d t 
= \e^{-r\eta_a(l+\eps)} \varsigma_a \eps . 
\label{eq: lem stag2}
\end{align}
By \eqref{eq: lem stag1}, we see that 
the function $ F(l) = \int_0^{\eta_a(l)} \e^{-rt} 1_{\{ X(t)=a \}} \d t $ 
is absolutely continuous, 
so that there exists a locally integrable function $ f(l) $ such that 
$ F(l) = \int_0^l f(s) \d s $. 
By \eqref{eq: lem stag1} and \eqref{eq: lem stag2}, we have 
the right-derivative of $ F(l) $ is equal to $ \varsigma_a \e^{- r\eta_a(l)} $. 
Hence we obtain $ f(l) = \varsigma_a \e^{- r\eta_a(l)} $ for almost every $ l \ge 0 $. 
The proof is now complete. 
}%endproof

For the proof of Theorem \ref{thm: main-1}, 
we prove three more lemmas in what follows. 

Let us construct a sample path of the process $ \{ (X_t)_{t \ge 0},P_a \} $. 
Let $ (p_a^{(b)}(l))_{l \in D_a^{(b)}} $ be a Poisson point process 
with characteristic measure $ n_a^{(b)} $ 
and let $ (X_b(t))_{t \ge 0} $ be a process with law $ P_b $. 
We assume that $ (p_a^{(b)}(l))_{l \ge 0} $ and $ (X_b(t))_{t \ge 0} $ are independent. 
We define the subordinator $ (\eta_a^{(b)}(l))_{l \ge 0} $ as 
\begin{align}
\eta_a^{(b)}(l) = \varsigma_a^{(b)} l + \sum_{s \le l} T_{a,b}(p_a^{(b)}(s)) 
\quad \text{for $ l \ge 0 $}. 
\label{}
\end{align}
Note that we have 
\begin{align}
E \sbra{ \e^{-r\eta_a^{(b)}(l)} } 
=& \exp \cbra{ - l \rbra{ \varsigma_a^{(b)} r + n^{(b)}_a \sbra{ 1-\e^{-r T_{a,b}} } } } 
\label{} \\
=& \exp \cbra{ - l \rbra{ \psi^{(b)}_a(r) - n^{(b)}_a \sbra{ \e^{-r T_b} ; T_b<\infty } } } . 
\label{}
\end{align}
Let $ \lambda $ denote the first index $ l \ge 0 $ 
that $ p_a^{(b)}(l) $ hits $ b $, or in other words, 
\begin{align}
\lambda = \inf \{ l \ge 0 : T_b(p_a^{(b)}(l)) < \infty \} . 
\label{}
\end{align}
We define the process $ (X_a^{(b)}(t) : 0 \le t < \eta_a^{(b)}(\lambda)) $ as 
\begin{align}
X_a^{(b)}(t) = 
\begin{cases}
p_a^{(b)}(l)(t-\eta_a^{(b)}(l-)) 
& \text{if $ \eta_a^{(b)}(l-) \le t < \eta_a^{(b)}(l) $ for some $ 0 \le l \le \lambda $}, \\
a & \text{otherwise}. 
\end{cases}
\label{}
\end{align}
Now we construct the process $ (X_a(t))_{t \ge 0} $ as 
\begin{align}
X_a(t) = 
\begin{cases}
X_a^{(b)}(t) 
& \text{if $ t < \eta_a^{(b)}(\lambda) $}, \\
X_b(t-\eta_a^{(b)}(\lambda)) 
& \text{if $ t \ge \eta_a^{(b)}(\lambda) $}. 
\end{cases}
\label{}
\end{align}
Then it is obvious that the process $ (X_a(t))_{t \ge 0} $ 
is a realization of the process $ \{ (X_t)_{t \ge 0},P_a \} $. 

The first one of the three lemmas is the following. 

\begin{Lem} \label{lem: main1}
Let $ g \in \cB_b $. 
Then 
\begin{align}
E \sbra{ \int_{\eta^{(b)}_a(\lambda)}^{\infty } \e^{-rt} g(X_a(t)) \d t } 
= \frac{n^{(b)}_a \sbra{ \e^{-r T_b} ; T_b < \infty }}{\psi_a^{(b)}(r)} 
\cdot R_rg(b) . 
\label{}
\end{align}
\end{Lem}

\Proof{
Let $ (p^0_a(l))_{l \in D_a^0} $ denote the restriction of $ (p^{(b)}_a(l))_{l \in D_a^{(b)}} $ 
on the set of excursions which fail to hit $ b $; more precisely, we define 
\begin{align}
D_a^0 = \{ l \in D_a^{(b)} : T_b(p^{(b)}_a(l)) = \infty \} 
\label{}
\end{align}
and $ p^0_a(l) = p^{(b)}_a(l) $ for all $ l \in D_a^0 $. 
Let $ \epsilon = p^{(b)}_a(\lambda) $. 
Then we see that the three quantities $ (p^0_a(l))_{l \ge 0} $, 
$ \lambda $ and $ \epsilon $ 
are mutually independent, 
and that 
$ (p^0_a(l))_{l \in D_a^0} $ is a Poisson point process 
with characteristic measure $ n^{(b)}_a(\cdot;T_b=\infty ) $, 
the law of $ \lambda $ is given as 
\begin{align}
P(\lambda >l) = \e^{- l \cdot n^{(b)}_a(T_b<\infty )} 
\label{}
\end{align}
and the law of $ \epsilon $ is given as 
\begin{align}
P(\epsilon \in \cdot) = \frac{n^{(b)}_a(\cdot;T_b<\infty )}{n^{(b)}_a(T_b<\infty )} . 
\label{}
\end{align}
If we write 
\begin{align}
\eta^0_a(l) = \varsigma_a^{(b)} l + \sum_{s \le l} T_a(p^0_a(s)) , 
\label{}
\end{align}
we have $ E[\e^{-r \eta^0_a(l)}] = \e^{-l \psi^0_a(r)} $, where 
\begin{align}
\psi^0_a(r) 
= \varsigma_a^{(b)} r + n^{(b)}_a \sbra{ 1-\e^{-r T_a} ; T_b = \infty } . 
\label{}
\end{align}
By definition, we have 
\begin{align}
\eta^{(b)}_a(\lambda) = \eta^0_a(\lambda) + T_b(\epsilon) . 
\label{eq: eta(1) and eta0}
\end{align}
Thus we obtain 
\begin{align}
E \sbra{ \int_{\eta^{(b)}_a(\lambda)}^{\infty } \e^{-rt} g(X_a(t)) \d t } 
=& 
E \sbra{ \e^{-r \eta^{(b)}_a(\lambda)} \int_0^{\infty } \e^{-rt} 
g \rbra{ X_a(t+\eta^{(b)}_a(\lambda)) \Big. } \d t } 
\label{} \\
=& 
E \sbra{ \e^{-r \eta^{(b)}_a(\lambda)} } 
E \sbra{ \int_0^{\infty } \e^{-rt} g(X_b(t)) \d t } 
\label{} \\
=& 
E \sbra{ \e^{-r \eta^0_a(\lambda) - r T_b(\epsilon)} } R_rg(b) 
\label{} \\
=& 
E \sbra{ \e^{-r \eta^0_a(\lambda) } } \cdot E \sbra{ \e^{- r T_b(\epsilon)} } R_rg(b) . 
\label{}
\end{align}
The expectations in the last expression can be computed as 
\begin{align}
E \sbra{ \e^{-r \eta^0_a(\lambda) } } 
= 
E \sbra{ \e^{ - \lambda \psi^0_a(r) } } 
= 
\frac{n^{(b)}_a(T_b < \infty )}{n^{(b)}_a(T_b < \infty ) + \psi^0_a(r)} 
\label{}
\end{align}
and 
\begin{align}
E \sbra{ \e^{- r T_b(\epsilon)} } 
= 
\frac{n^{(b)}_a \sbra{ \e^{-r T_b} ; T_b < \infty }}{n^{(b)}_a(T_b < \infty )} . 
\label{}
\end{align}
Since 
\begin{align}
n^{(b)}_a(T_b < \infty ) + \psi^0_a(r) 
= \varsigma_a^{(b)} r + n^{(b)}_a \sbra{ 1-\e^{-r T_a} } 
= \psi_a^{(b)}(r) , 
\label{}
\end{align}
we complete the proof. 
}

The second one is the following. 

\begin{Lem} \label{lem: main2}
Let $ g \in \cB_b $. Then 
\begin{align}
E \sbra{ \int_0^{\eta^{(b)}_a(\lambda)} \e^{-rt} g(X_a(t)) 1_{\{ X_a(t)=a \}} \d t } 
= \varsigma_a^{(b)} g(a) \frac{1}{\psi_a^{(b)}(r)} . 
\label{eq: lem main2-1}
\end{align}
\end{Lem}

\Proof{
By Lemma \ref{lem: stagnancy}, we have 
\begin{align}
\text{the left-hand side of \eqref{eq: lem main2-1}} 
=& g(a) 
E \sbra{ \int_0^{\eta^{(b)}_a(\lambda)} \e^{-rt} 1_{\{ X_a(t)=a \}} \d t } 
\label{} \\
=& g(a) 
E \sbra{ \varsigma_a^{(b)} \int_0^{\lambda} \e^{-r \eta^{(b)}_a(l)} \d l } 
\label{} \\
=& \varsigma_a^{(b)} g(a) 
E \sbra{ \int_0^{\lambda} \e^{-r \eta^0_a(l)} \d l } . 
\label{}
\end{align}
Since $ \lambda $ is independent of $ (\eta^0_a(l))_{l \ge 0} $, we have 
\begin{align}
E \sbra{ \int_0^{\lambda} \e^{-r \eta^0_a(l)} \d l } 
=& 
E \sbra{ \int_0^{\lambda} \e^{-l \psi^0_a(r)} \d l } 
\label{eq: lem main2-2} \\
=& \frac{1}{\psi^0_a(r)} 
E \sbra{ 1 - \e^{- \lambda \psi^0_a(r)} } 
\label{} \\
=& \frac{1}{\psi^0_a(r)} 
\cbra{ 1 - \frac{n^{(b)}_a(T_b<\infty )}{n^{(b)}_a(T_b<\infty ) + \psi^0_a(r)} } 
\label{} \\
=& \frac{1}{n^{(b)}_a(T_b<\infty ) + \psi^0_a(r)} 
\label{} \\
=& \frac{1}{\psi^{(b)}_a(r)} , 
\label{eq: lem main2-3}
\end{align}
which completes the proof. 
}

The third one is the following. 

\begin{Lem} \label{lem: main3}
Let $ g \in \cB_b $. Then 
\begin{align}
E \sbra{ \int_0^{\eta^{(b)}_a(\lambda)} \e^{-rt} g(X_a(t)) \d t } 
= \frac{1}{\psi_a^{(b)}(r)} N^{(b)}_{a,r}(g) . 
\label{eq: lem main3-1}
\end{align}
\end{Lem}

\Proof{
By the construction of the process $ X_a(t) $, we have 
\begin{align}
& E \sbra{ \int_0^{\eta^{(b)}_a(\lambda)} \e^{-rt} g(X_a(t)) 1_{\{ X_a(t) \neq a \}} \d t } 
\label{} \\
=& E \sbra{ \sum_{l \le \lambda} 
\int_{\eta^{(b)}_a(l-)}^{\eta^{(b)}_a(l)} \e^{-rt} g(X_a(t)) \d t } 
\label{} \\
=& E \sbra{ \sum_{l \le \lambda} \e^{-r \eta^{(b)}_a(l-)} 
\int_{0}^{T_{a,b}(p^{(b)}_a(l))} \e^{-rt} g(p^{(b)}_a(l)(t)) \d t } . 
\label{eq: lem main3-2}
\end{align}
By the compensation formula, we have 
\begin{align}
\eqref{eq: lem main3-2} 
=& E \sbra{ \int_0^{\lambda} \e^{-r \eta^{(b)}_a(l)} \d l } 
\cdot n^{(b)}_a \sbra{ \int_0^{T_{a,b}} \e^{-rt} g(X_t) \d t } . 
\label{}
\end{align}
By \eqref{eq: eta(1) and eta0} and by \eqref{eq: lem main2-2}-\eqref{eq: lem main2-3}, we have 
\begin{align}
E \sbra{ \int_0^{\lambda} \e^{-r \eta^{(b)}_a(l)} \d l } 
= E \sbra{ \int_0^{\lambda} \e^{-r \eta^0_a(l)} \d l } 
= \frac{1}{\psi_a^{(b)}(r)} . 
\label{}
\end{align}
Combining these results with the result of Lemma \ref{lem: main2}, 
we obtain the desired result. 
}

Theorem \ref{thm: main-1} is therefore immediate 
from Lemmas \ref{lem: main1} and \ref{lem: main3}.

\section{The resolvent of the minimal diffusion} \label{sec: minimal}

\subsection{Non-negative increasing and decreasing eigenfunctions}

Let $ \{ (X_t)_{t \ge 0},(P^0_x)_{x \in (a,b)} \} $ 
be a minimal diffusion process 
and follow the notations in Subsection \ref{sec: mp}. 

We recall non-negative increasing and decreasing eigenfunctions of $ \cL $. 
All results in this subsection are well-known; 
see, for example, \cite{IM} for details. 

Let $ c \in (a,b) $ and $ r>0 $ be fixed. 
Let $ v=\varphi_{r} $ and $ \psi_{r} $ denote 
the unique non-negative increasing solutions of $ \cL v = r v $ such that 
\begin{align}
& \varphi_{r}(c) = 1, 
\quad 
\cD_s \varphi_{r}(c) = 0 , 
\label{} \\
& \psi_{r}(c) = 0 , 
\quad 
\cD_s \psi_{r}(c) = 1 . 
\label{}
\end{align}
These solutions can be obtained via successive approximation 
by solving the following integral equations: 
\begin{align}
\varphi_{r}(x) =& 1 + r \int_c^x \d s(y) \int_c^y \varphi_{r}(z) \d m(z) 
, \label{} \\
\psi_{r}(x) =& s(x)-s(c) + r \int_c^x \d s(y) \int_c^y \psi_{r}(z) \d m(z) . 
\label{}
\end{align}
Then any solution of $ \cL v = r v $ with $ v(c)=1 $ 
is of the form $ v = \varphi_{r} - \gamma \psi_{r} $. 
Note that $ v = \varphi_{r} - \gamma \psi_{r} $ 
is non-negative and decreasing when restricted on $ [c,b) $, 
then so is it on the whole interval $ (a,b) $. 
Let $ \underline{\gamma} $ and $ \bar{\gamma} $ denote 
the infimum and the supremum among all $ \gamma>0 $ such that 
$ \varphi_{r} - \gamma \psi_{r} $ is a non-negative decreasing function. 
Then we see that $ 0 < \underline{\gamma} \le \bar{\gamma} < \infty $ 
and that $ \varphi_{r} - \gamma \psi_{r} $ is a non-negative decreasing function 
for all $ \underline{\gamma} \le \gamma \le \bar{\gamma} $. 
Now we take the minimal one: $ v_{r} = \varphi_{r} - \bar{\gamma} \psi_{r} $. 
The boundary behaviors of $ v_{r} $ at $ a $ are as follows: 
\begin{center}
\begin{tabular}{l|l}
$ v_{r}(a) \begin{cases}
            \in (0,\infty ) & \text{if $ a $ is accessible} \\
            = \infty & \text{otherwise} 
           \end{cases}$
                                           & $ v_{r}(b) 
                                              \begin{cases}
                                               \in (0,\infty ) & \text{if $ b $ is entrance} \\
                                               0 & \text{otherwise}
                                              \end{cases}
                                             $   
\\ \hline
$ -\cD_s v_{r}(a) \begin{cases}
 \in (0,\infty ) & \text{if $ a $ is enterable} \\
 = \infty & \text{otherwise} 
 \end{cases}$
                                           & $ -\cD_s v_{r}(b) 
                                              \begin{cases}
                                               \in (0,\infty ) & \text{if $ b $ is accessible} \\
                                               = 0 & \text{otherwise} 
                                              \end{cases}
                                             $ 
\end{tabular}
\end{center}
We also note that $ \int_c^{b-} v_r(x) \d m(x) < \infty $ in any case. 
By the same way, we obtain a non-negative increasing solution 
$ u=u_{r} $ of $ \cL u = r u $ 
whose boundary behaviors are as follows: 
\begin{center}
\begin{tabular}{l|l}
$ u_{r}(a) \begin{cases}
            \in (0,\infty ) & \text{if $ a $ is entrance} \\
            0 & \text{otherwise}
           \end{cases} $
                            & $ u_{r}(b) \begin{cases}
                                          \in (0,\infty ) & \text{if $ b $ is accessible} \\
                                          = \infty & \text{otherwise} 
                                         \end{cases}$
\\ \hline
$ \cD_s u_{r}(a) \begin{cases}
                \in (0,\infty ) & \text{if $ a $ is accessible} \\
                = 0 & \text{otherwise} 
               \end{cases}$
                        & $ \cD_s u_{r}(b) \begin{cases}
                                          \in (0,\infty ) & \text{if $ b $ is enterable} \\
                                          = \infty & \text{otherwise} 
                                         \end{cases}$
\end{tabular}
\end{center}
We also note that $ \int_a^c u_r(x) \d m(x) < \infty $ in any case.

\subsection{Several limits at the boundary points} \label{sec: bdry limit}

We multiply $ u_r $ (or $ v_r $) by a certain constant and 
we may assume without loss of generality that 
\begin{align}
v_r(x) \cD_su_r(x) - u_r(x) \cD_sv_r(x) = 1 
\quad \text{for all $ x \in (a,b) $}. 
\label{eq: W}
\end{align}
We prove the following proposition. 

\begin{Prop} \label{prop: BL}
One has 
\begin{align}
\cD_su_r(a) = \frac{1}{v_r(a)} 
\quad \text{and} \quad 
\cD_sv_r(b) = - \frac{1}{u_r(b)} . 
\label{}
\end{align}
Here we understand $ \frac{1}{+\infty } = 0 $ and $ \frac{1}{+0} = +\infty $. 
\end{Prop}

\Proof{
By the symmetry, it suffices only to prove that $ \cD_su_r(a) = 1/v_r(a) $. 
If $ a $ is not accessible, 
then this is obvious, 
because we have $ \cD_su_r(a) = 0 $ and $ v_r(a) = \infty $. 

Suppose that $ a $ is accessible. 
By \eqref{eq: W}, we have 
\begin{align}
\cD_su_r(a) = \lim_{x \to a+} \frac{1 - u_r(x) \cD_sv_r(x)}{v_r(x)} . 
\label{}
\end{align}
Hence it suffices to show that 
\begin{align}
\lim_{x \to a+} \{ - u_r(x) \cD_sv_r(x) \} = 0 . 
\label{eq: lim of urDvr}
\end{align}
Since $ u_r $ satisfies $ \cL u_r = r u_r $ and $ u_r(a)=0 $, we have 
\begin{align}
u_r(x) = k s(a,x) + r \int_a^x \d s(y) \int_a^y u_r(z) \d m(z) , 
\label{}
\end{align}
where we denote $ k = \cD_su_r(a) \in (0,\infty ) $. 
Differentiating both sides, we have 
\begin{align}
\cD_su_r(x) = k + r \int_a^x u_r(z) \d m(z) . 
\label{}
\end{align}
Let $ \eps>0 $ be fixed. 
Then there exists $ \delta>0 $ such that 
$ |\cD_su_r(x) - k| \le \eps $ for all $ x $ with $ a<x<a+\delta $. 
Then we have 
\begin{align}
|u_r(x) - ks(a,x)| \le \int_a^x |\cD_su_r(y)-k| \d s(y) \le \eps s(a,x) 
\quad \text{for} \ a<x<a+\delta. 
\label{}
\end{align}
By \eqref{eq: W}, we see that $ \cD_s(v_r/u_r) = -1/u_r^2 $ and that 
\begin{align}
\frac{v_r(x)}{u_r(x)} 
= \frac{v_r(b)}{u_r(b)} + \int_x^{b-} \frac{\d s(y)}{u_r(y)^2} 
= \int_x^{b-} \frac{\d s(y)}{u_r(y)^2} . 
\label{}
\end{align}
From this and by \eqref{eq: W}, we have 
\begin{align}
0 \le - u_r(x) \cD_sv_r(x) 
=& 1 - v_r(x) \cD_s u_r(x) 
\label{} \\
=& 1 - u_r(x) \cD_su_r(x) \int_x^{b-} \frac{\d s(y)}{u_r(y)^2} 
\label{} \\
\le& 1 - (k-\eps)^2 s(a,x) \int_x^{b-} \frac{\d s(y)}{\{ (k+\eps) s(a,y) \}^2} 
\label{} \\
\le& 1 - \frac{(k-\eps)^2}{(k+\eps)^2} s(a,x) \cbra{ \frac{1}{s(a,x)} - \frac{1}{s(a,b)} } 
\label{} \\
\le& 1 - \frac{(k-\eps)^2}{(k+\eps)^2} 
+ \frac{(k-\eps)^2}{(k+\eps)^2} \cdot \frac{s(a,x)}{s(a,b)} . 
\label{}
\end{align}
Since $ s(a,a)=0 $ and since $ \eps>0 $ is arbitrary, 
we obtain \eqref{eq: lim of urDvr} 
and hence we obtain the desired result. 
}

\subsection{The resolvent of the minimal diffusion process}

Define 
\begin{align}
R^0_r(x,y) = R^0_r(y,x) = u_{r}(x) v_{r}(y) 
\quad \text{for} \ a < x \le y < b . 
\label{}
\end{align}
Then it is well-known that 
the resolvent operator $ (R^0_r)_{r>0} $ of the minimal diffusion process 
defined in \eqref{eq: resolv of mdp} 
is the integral operator with kernel $ R^0_r(x,y) \d m(y) $, i.e., 
\begin{align}
R^0_rg(x) = \int_a^{b-} R^0_r(x,y) g(y) \d m(y) 
\quad \text{for $ g \in \cB_b $ and $ x \in (a,b) $}. 
\label{}
\end{align}
From Theorem \ref{thm: g of mp}, it follows that 
\begin{align}
\text{if $ f \in D(\cG^0) $, then }
R^0_r \cL f = r R^0_r f - f . 
\label{eq: RL0}
\end{align}
In our study, however, we need $ R^0_r \cL f $ for $ f \in D(\cL) $ with $ f,\cL f \in C_b $. 
The formula \eqref{eq: RL0} can be generalized to the following proposition: 

\begin{Prop} \label{prop: RL}
For any $ f \in D(\cL) $ with $ f,\cL f \in C_b $, one has 
\begin{align}
R^0_r \cL f = r R^0_r f - f + f(a) \frac{v_r}{v_r(a)} + f(b) \frac{u_r}{u_r(b)} . 
\label{eq: RL formula}
\end{align}
In addition, one has 
\begin{align}
\cD_sf(a)=0 \ \text{if $ a $ is entrance}. 
\label{}
\end{align}
\end{Prop}

\Proof{
Since $ \cL = \cD_m \cD_s $, we have 
\begin{align}
\cD_sf(x)-\cD_sf(y) = \int_y^x \cL f(z) \d m(z) 
\quad \text{for} \ x,y \in (a,b) . 
\label{eq: f'}
\end{align}
Since $ \cL u_r = r u_r $ and $ \cL v_r = r v_r $, we see that, 
for any $ x,y \in (a,b) $, 
\begin{align}
\cD_su_r(x) - \cD_su_r(y) = r \int_y^x u_r(z) \d m(z) , 
\label{eq: u'} \\
\cD_sv_r(x) - \cD_sv_r(y) = r \int_y^x v_r(z) \d m(z) . 
\label{eq: v'}
\end{align}

Let $ x \in (a,b) $ 
and take $ c \in (a,x) $ arbitrarily. 
We have 
\begin{align}
& \int_c^x \cL f(y) u_r(y) \d m(y) 
\label{} \\
=& \int_c^x \d m(y) \cL f(y) \cbra{ u_r(c) + \int_c^y \cD_su_r(z) \d s(z) } 
\label{} \\
=& u_r(c) \cbra{ \cD_sf(x)-\cD_sf(c) } + \int_c^x \d s(z) \cD_su_r(z) \int_z^x \cL f(y) \d m(y) 
\label{} \\
=& u_r(c) \cbra{ \cD_sf(x)-\cD_sf(c) } + \int_c^x \d s(z) \cD_su_r(z) \{ \cD_sf(x)-\cD_sf(z) \} 
\label{} \\
=& u_r(x) \cD_sf(x) - u_r(c) \cD_sf(c) - \int_c^x \d s(z) \cD_su_r(z) \cD_sf(z) . 
\label{}
\end{align}
Using the formula \eqref{eq: u'}, we have 
\begin{align}
& \int_c^x \d s(z) \cD_su_r(z) \cD_sf(z)
\label{} \\
=& \int_c^x \d s(z) \cD_sf(z) \{ \cD_su_r(x) - (\cD_su_r(x) - \cD_su_r(z)) \} 
\label{} \\
=& \cD_su_r(x) \{ f(x)-f(c) \} - r \int_c^x \d s(z) \cD_sf(z) \int_z^x u_r(y) \d m(y) 
\label{} \\
=& \cD_su_r(x) \{ f(x)-f(c) \} - r \int_c^x \d m(y) u_r(y) \int_c^y \cD_sf(z) \d s(z) 
\label{} \\
=& \cD_su_r(x) \{ f(x)-f(c) \} - r \int_c^x f(y) u_r(y) \d m(y) 
+ f(c) \cdot r \int_c^x u_r(y) \d m(y) . 
\label{eq: RL1}
\end{align}
Using the formula \eqref{eq: u'} again, we have 
\begin{align}
\text{\eqref{eq: RL1}}
=& \cD_su_r(x) \{ f(x)-f(c) \} - r \int_c^x f(y) u_r(y) \d m(y) 
+ f(c) \{ \cD_su_r(x)-\cD_su_r(c) \} 
\label{} \\
=& f(x) \cD_su_r(x) - f(c) \cD_su_r(c) - r \int_c^x f(y) u_r(y) \d m(y) . 
\label{}
\end{align}
Hence we obtain 
\begin{align}
\int_c^x \cL f(y) u_r(y) \d m(y) 
= W[u_r,f](x) - W[u_r,f](c) + r \int_c^x f(y) u_r(y) \d m(y) , 
\label{eq: RL2}
\end{align}
where we write 
\begin{align}
W[f,g](x) = f(x) \cD_s g(x) - g(x) \cD_s f(x) . 
\label{}
\end{align}
Since $ f $ and $ \cL f $ are bounded and since $ \int_a^c u_r(y) \d m(y) $ is finite, 
we see that the limit 
\begin{align}
Q_a := \lim_{z \to a+} \cbra{ - W[u_r,f](z) } 
\label{}
\end{align}
exists finitely and that 
\begin{align}
\int_a^x \cL f(y) u_r(y) \d m(y) 
= Q_a + W[u_r,f](x) + r \int_a^x f(y) u_r(y) \d m(y) 
\label{eq: RL3}
\end{align}
holds. In the same way, we see that the limit 
\begin{align}
Q_b := \lim_{z \to b-} W[v_r,f](z) 
\label{}
\end{align}
exists finitely and that 
\begin{align}
\int_x^{b-} \cL f(y) v_r(y) \d m(y) 
= Q_b - W[v_r,f](x) + r \int_x^{b-} f(y) v_r(y) \d m(y) 
\label{eq: RL4}
\end{align}
holds. 
Adding \eqref{eq: RL3} times $ v_r(x) $ and \eqref{eq: RL4} times $ u_r(x) $, 
we obtain 
\begin{align}
(R^0_r \cL f)(x) = Q_a v_r(x) + Q_b u_r(x) - f + r R^0_rf(x) . 
\label{eq: RL5}
\end{align}

Now let us prove that $ Q_a = f(a)/v_r(a) $ as follows: 
\begin{enumerate}
\item 
Suppose that $ a $ is accessible. Let $ \| g \| = \sup_x |g(x)| < \infty $. 
Since $ v_r $ is decreasing and by \eqref{eq: lim of urDvr}, 
we have 
\begin{align}
u_r(c) | \cD_sf(x)-\cD_sf(c) | 
\le& u_r(c) \int_c^x |\cL f(z)| \d m(z) 
\label{} \\
\le& \frac{\| \cL f \|}{v_r(x)} u_r(c) \int_c^x v_r(z) \d m(z) 
\label{} \\
=& \frac{\| \cL f \|}{r v_r(x)} u_r(c) \{ \cD_s v_r(x) - \cD_s v_r(c) \} 
\label{} \\
\to& 0 \quad \text{as $ c \to a+ $ for fixed $ x $}. 
\label{}
\end{align}
Hence we obtain $ \lim_{c \to a+} \{ - u_r(c) \cD_s f(c) \} = 0 $. 
Hence, by Proposition \ref{prop: BL}, we obtain 
$ Q_a = f(a) \cD_su_r(a) = f(a)/v_r(a) $. 
\item 
Suppose that $ a $ is not accessible. 
On one hand, we have $ v_r(a) = \infty $. 
On the other hand, since $ f $ and $ \cL f $ are bounded, 
we see by \eqref{eq: RL5} that $ Q_a v_r(x) $ should be bounded near $ a $. 
Hence we obtain $ Q_a = 0 = f(a)/v_r(a) $. 
\end{enumerate}
We can make the same argument for $ b $ and obtain $ Q_b = f(b)/u_r(b) $. 
Therefore, from \eqref{eq: RL5}, we obtain the formula \eqref{eq: RL formula}. 

If $ a $ is entrance, 
then we have $ u_r(a) \in (0,\infty ) $ and $ \cD_s u_r(a) = 0 $. 
Since $ Q_a=0 $, we obtain $ \cD_s f(a) = 0 $. 

The proof is now complete. 
}

\section{The $ C_b $-generator} \label{sec: determ}

We now suppose that 
$ \{ (X_t)_{t \ge 0},(P_x)_{x \in [a,b]} \} $ is 
an extension to a standard process of the minimal process 
$ \{ (X_t)_{t \ge 0},(P^0_x)_{x \in (a,b)} \} $ 
and follow the notations in Subsection \ref{sec: emp}. 
It is well-known that 
\begin{align}
E^{\rm stop}_x \sbra{ \e^{-rT_a} ; T_a<T_b } = \frac{v_r(x)}{v_r(a)} 
\quad \text{and} \quad 
E^{\rm stop}_x \sbra{ \e^{-rT_b} ; T_a>T_b } = \frac{u_r(x)}{u_r(b)} 
\quad \text{for} \ r>0 . 
\label{eq: stop hitting}
\end{align}

\subsection{$ \Phi_a(v_r) $ and $ \Phi_a(u_r) $}

We need the following lemma for later use. 

\begin{Lem} \label{lem: Phiavr} 
If $ a $ is accessible, one has 
\begin{align}
\frac{\Phi_a(v_r)}{v_r(a)} 
=& p_3 r + n^{(b)}_a \sbra{ 1 - \e^{-r T_a} } = \psi^{(b)}_a(r) . 
\label{eq: Phiavr1}
\end{align}
If $ a $ and $ b $ are both accessible, one has 
\begin{align}
\frac{\Phi_a(u_r)}{u_r(b)} 
=& - n^{(b)}_a \sbra{ \e^{-r T_b} ; T_b < \infty } . 
\label{eq: Phiavr2}
\end{align}
\end{Lem}

\Proof{
By \eqref{eq: nrefl} and \eqref{eq: stop hitting}, we have 
\begin{align}
- \frac{\cD_s v_r(a)}{v_r(a)} 
=& \lim_{x \to a+} \frac{1}{s(a,x)} \cbra{ 1-\frac{v_r(x)}{v_r(a)} } 
\label{} \\
=& \lim_{x \to a+} n^{(b)}_{a,{\rm refl}} \rbra{ T_x < \infty } 
\cbra{ 1 - P^{\rm stop}_x \sbra{ \e^{- r T_a} ; T_a < \infty } } 
\label{} \\
=& \lim_{x \to a+} n^{(b)}_{a,{\rm refl}} \rbra{ T_x < \infty } 
\cbra{ P^{\rm stop}_x \sbra{ 1-\e^{- r T_a} ; T_a < \infty } + P^{\rm stop}_x(T_b < \infty ) } . 
\label{eq: prf of Lem1}
\end{align}
Note that, under the measure $ n^{(b)}_{a,{\rm refl}} $, 
the hitting time $ T_x $ decreases to 0 as $ x $ does to $ a $. 
By the strong Markov property of $ n^{(b)}_{a,{\rm refl}} $ 
and by the dominated convergence theorem, 
we have 
\begin{align}
\eqref{eq: prf of Lem1} 
=& \lim_{x \to a+} \cbra{ 
n^{(b)}_{a,{\rm refl}} \sbra{ 1-\e^{- r (T_a-T_x)} ; T_x < T_a < \infty } 
+ n^{(b)}_{a,{\rm refl}} \rbra{ T_x < T_b < \infty } } 
\label{} \\
=& n^{(b)}_{a,{\rm refl}} \sbra{ 1-\e^{- r T_a};T_a < \infty } 
+ n^{(b)}_{a,{\rm refl}}(T_b < \infty ) 
\label{} \\
=& n^{(b)}_{a,{\rm refl}} \sbra{ 1-\e^{- r T_a} } . 
\label{}
\end{align}
Hence, by \eqref{eq: n(b)a}, we have 
\begin{align}
- p_2 \frac{\cD_s v_r(a)}{v_r(a)} 
= n^{(b)}_a \sbra{ 1-\e^{- r T_a} ; X_0=a } . 
\label{}
\end{align}
By \eqref{eq: stop hitting} and by \eqref{eq: n(b)a}, we have 
\begin{align}
& \int_{(a,b]} p_4(\d x) \cbra{ 1-\frac{v_r(x)}{v_r(a)} } 
\label{} \\
=& 
\int_{(a,b]} p_4(\d x) \cbra{ P^{\rm stop}_x \sbra{ 1-\e^{-r T_a} ; T_a < \infty } 
+ P^{\rm stop}_x(T_b < \infty ) } 
\label{} \\
=& n^{(b)}_a \sbra{ 1-\e^{-r T_a} ; T_a < \infty , \ X_0 \in (a,b] } 
+ n^{(b)}_a (T_b < \infty , \ X_0 \in (a,b] ) . 
\label{}
\end{align}
Since $ p_1 = n_a(\{ \Delta \}) $ 
and since $ \{ \Delta \} \cup \{ T_b < \infty \} = \{ T_a=\infty \} $, we obtain 
\begin{align}
\frac{\Phi_a(v_r)}{v_r(a)} 
=& p_1 - p_2 \frac{\cD_s v_r(a)}{v_r(a)} + p_3 \frac{\cL v_r(a)}{v_r(a)} 
+ \int_{(a,b]} p_4(\d x) \cbra{ 1-\frac{v_r(x)}{v_r(a)} } 
\label{} \\
=& p_3 r + n^{(b)}_a \sbra{ 1 - \e^{-r T_a} ; T_a < \infty } 
+ n^{(b)}_a(\{ \Delta \}) + n^{(b)}_a \rbra{ T_b < \infty } 
\label{} \\
=& p_3 r + n^{(b)}_a \sbra{ 1 - \e^{-r T_a} } . 
\label{}
\end{align}
Now we obtain \eqref{eq: Phiavr1}. 

Since $ a $ is regular, we have $ u_r(a)=0 $. 
By \eqref{eq: nrefl} and \eqref{eq: stop hitting}, we have 
\begin{align}
\frac{\cD_s u_r(a)}{u_r(b)} 
=& \lim_{x \to a+} \frac{1}{s(a,x)} \frac{u_r(x)}{u_r(b)} 
\label{} \\
=& \lim_{x \to a+} n^{(b)}_{a,{\rm refl}} \rbra{ T_x < \infty } 
P^{\rm stop}_x \sbra{ \e^{-r T_b} ; T_b < \infty } 
\label{} \\
=& \lim_{x \to a+} n^{(b)}_{a,{\rm refl}} \sbra{ \e^{-r(T_b-T_x)} ; T_x < T_b < \infty } 
\label{} \\
=& n^{(b)}_{a,{\rm refl}} \sbra{ \e^{-r T_b} ; T_b < \infty } . 
\label{}
\end{align}
Hence we have 
\begin{align}
p_2 \frac{\cD_s u_r(a)}{u_r(b)} 
= n^{(b)}_a \sbra{ \e^{-r T_b} ; T_b < \infty , \ X_0=a } . 
\label{}
\end{align}
By \eqref{eq: stop hitting} and by \eqref{eq: n(b)a}, we have 
\begin{align}
\int_{(a,b]} p_4(\d x) \frac{u_r(x)}{u_r(a)} 
=& 
\int_{(a,b]} p_4(\d x) P^{\rm stop}_x \sbra{ \e^{-r T_b} ; T_b < \infty } 
\label{} \\
=& n^{(b)}_a \sbra{ \e^{-r T_b} ; \ T_b < \infty , \ X_0 \in (a,b] } . 
\label{}
\end{align}
Since $ u_r(a) = 0 $, we have 
\begin{align}
\frac{\Phi_a(u_r)}{u_r(b)} 
=& p_1 \frac{u_r(a)}{u_r(b)} - p_2 \frac{\cD_s u_r(a)}{u_r(b)} + p_3 \frac{\cL u_r(a)}{u_r(b)} 
- \int_{(a,b]} p_4(\d x) \frac{u_r(x)-u_r(a)}{u_r(b)} 
\label{} \\
=& - n^{(b)}_a \sbra{ \e^{-r T_b} ; T_b < \infty } . 
\label{}
\end{align}
The proof is now complete. 
}

The following proposition is important 
in the proof of our main theorem. 

\begin{Prop} \label{prop: main2}
If $ a $ and $ b $ are both accessible, one has 
\begin{align}
\frac{\Phi_a(v_r)}{v_r(a)} \cdot \frac{\Phi_b(u_r)}{u_r(b)} 
- 
\frac{\Phi_b(v_r)}{v_r(a)} \cdot \frac{\Phi_a(u_r)}{u_r(b)} 
> 0 
\quad \text{for} \ r>0 . 
\label{}
\end{align}
That is, the matrix 
\begin{align}
A = \pmat{ \Phi_a(v_r)/v_r(a) & \Phi_a(u_r)/u_r(b) \\ \Phi_b(v_r)/v_r(a) & \Phi_b(u_r)/u_r(b) } 
\label{eq: A matrix}
\end{align}
has strictly positive determinant. 
\end{Prop}

\Proof{
Let us write $ F_x = (1-\e^{-r T_x}) 1_{\{ T_x<\infty \}} $ for $ x=a $ and $ b $. 
By Lemma \ref{lem: Phiavr} and by that where $ a $ and $ b $ are switched, we have 
\begin{align}
& \frac{\Phi_a(v_r)}{v_r(a)} 
\cdot \frac{\Phi_b(u_r)}{u_r(b)} 
- 
\frac{\Phi_b(v_r)}{v_r(a)} 
\cdot \frac{\Phi_a(u_r)}{u_r(b)} 
\label{} \\
\ge& 
\cbra{ p_3 r + n^{(b)}_a[F_a] + n^{(b)}_a(\{ \Delta \}) \Big. } 
\cdot 
\cbra{ q_3 r + n^{(a)}_b[F_b] + n^{(a)}_b(\{ \Delta \}) } 
\n \\
&+ 
n^{(b)}_a(T_b < \infty ) \cdot n^{(a)}_b(T_a < \infty ) 
- n^{(b)}_a \sbra{ \e^{-r T_b} ; T_b < \infty } 
\cdot n^{(a)}_b \sbra{ \e^{-r T_b} ; T_a < \infty } 
\label{} \\
\ge& 
\cbra{ p_3 r + p_2 n^{(b)}_{a,{\rm refl}}[F_a] + p_1 
+ \int_{(a,b)} p_4(\d x) P^{\rm stop}_x[F_a] } 
\n \\
&\cdot 
\cbra{ q_3 r + q_2 n^{(a)}_{b,{\rm refl}}[F_b] + q_1 
+ \int_{(a,b)} q_4(\d x) P^{\rm stop}_x[F_b] } . 
\label{}
\end{align}
The last quantity turns out to be positive 
because of the conditions \eqref{eq: pcond2} and \eqref{eq: qcond2}. 
}

\subsection{$ N^{(b)}_{a,r}(g) $ and $ \Phi_a(R^0_rg) $}

Let us prove the following lemma. 

\begin{Lem} \label{lem: main4}
Suppose that $ a $ is accessible. 
Then, for any $ g \in \cB_b $, one has 
\begin{align}
N^{(b)}_{a,r}(g) = - \Phi_a(R^0_rg) . 
\label{}
\end{align}
\end{Lem}

\Proof{
By Proposition \ref{prop: L0R}, we note that $ R^0_rg(a)=0 $. 

Suppose $ a $ is regular for a while. 
By the strong Markov property of $ n^{(b)}_{a,{\rm refl}} $, we have 
\begin{align}
\cD_s R^0_rg(a) 
=& \lim_{x \to a+} \frac{1}{s(a,x)} R^0_rg(x) 
\label{} \\
=& \lim_{x \to a+} n^{(b)}_{a,{\rm refl}}(T_x < \infty ) 
P^{\rm stop}_x \sbra{ \int_0^{T_{a,b}} \e^{-rt} g(X_t) \d t } 
\label{} \\
=& \lim_{x \to a+} n^{(b)}_{a,{\rm refl}} \sbra{ \int_0^{T_{a,b}-T_x} \e^{-rt} g(X_{t+T_x}) \d t 
\ ; \ T_x < \infty } 
\label{} \\
=& \lim_{x \to a+} n^{(b)}_{a,{\rm refl}} \sbra{ \int_{T_x}^{T_{a,b}} \e^{-r(t-T_x)} g(X_t) \d t 
\ ; \ T_x < \infty } . 
\label{}
\end{align}
Since 
\begin{align}
\absol{ \int_{T_x}^{T_{a,b}} \e^{-r(t-T_x)} g(X_t) \d t } 
\le \frac{\| g \|}{r} \rbra{ 1-\e^{-r T_{a,b}} } 
\label{}
\end{align}
and since $ n^{(b)}_{a,{\rm refl}}[1-\e^{-r T_{a,b}}] < \infty $, 
we may apply the dominated convergence theorem to see that 
\begin{align}
\cD_s R^0_rg(a) 
= n^{(b)}_{a,{\rm refl}} \sbra{ \int_0^{T_{a,b}} \e^{-rt} g(X_t) \d t } . 
\label{}
\end{align}

By Proposition \ref{prop: L0R}, we have 
\begin{align}
\cL R^0_rg(a) = r R^0_rg(a) - g(a) = -g(a) . 
\label{}
\end{align}
Therefore we obtain 
\begin{align}
- \Phi_a(R^0_rg) 
=& - p_1 R^0_rg(a) + p_2 \cD_sR^0_rg(a) - p_3 \cL R^0_rg(a) 
+ p_4[R^0_rg-R^0_rg(a)] 
\label{} \\
=& p_3 g(a) 
+ \rbra{ p_2 n^{\rm refl}_a + \int_{(a,b]} p_4(\d x) P^{\rm stop}_x } 
\sbra{ \int_0^{T_{a,b}} \e^{-rt} g(X_t) \d t } 
\label{} \\
=& p_3 g(a) + n^{(b)}_a \sbra{ \int_0^{T_{a,b}} \e^{-rt} g(X_t) \d t } . 
\label{}
\end{align}
The proof is now complete. 
}

\subsection{The case of {\boldmath $ 1^{\circ}). $}}

Let us prove Theorem \ref{thm: main--}

\Proof[Proof of Theorem \ref{thm: main--}]{
Suppose that $ a $ is accessible and that $ b $ is not. 
Let $ g \in \cB_b = \cB_b([a,b)) $. 
Noting that the process cannot hit $ b $ before hitting $ a $, we have the Dynkin formula: 
\begin{align}
R_rg = R^0_rg + R_rg(a) \frac{v_r}{v_r(a)} 
\quad \text{on $ [a,b) $}. 
\label{eq: Dynkin--}
\end{align}
By Proposition \ref{prop: L0R}, we have $ R_rg \in D(\cL) $, and we have 
\begin{align}
\cL R_rg 
= rR^0_rg-g + R_rg(a) \frac{rv_r}{v_r(a)} 
= r R_rg - g 
\quad \text{on $ [a,b) $}. 
\label{eq: rRrg-g--}
\end{align}
Using Lemmas \ref{lem: Phiavr} and \ref{lem: main4} 
and then using Theorem \ref{thm: main-1}, we have 
\begin{align}
\Phi_a(R_rg) 
=& \Phi_a(R^0_rg) + R_rg(a) \frac{\Phi_a(v_r)}{v_r(a)} 
\label{eq: main--1} \\
=& - N^{(b)}_{a,r}(g) + R_rg(a) \psi^{(b)}_a(r) 
\label{} \\
=& R_rg(b) n^{(b)}_a[\e^{-rT_b};T_b<\infty ] 
\label{eq: main--2} \\
=& 0 . 
\label{}
\end{align}
Thus we obtain the following: 
\begin{align}
\text{if $ g \in \cB_b $, we have $ R_rg \in D(\cL) $ and $ \Phi_a(R_rg)=0 $}. 
\label{eq: LR1--}
\end{align}
Set 
\begin{align}
\tilde{D} = \cbra{ f \in D(\cL) : 
f, \cL f \in C_b([a,b)) \ \text{and} \ 
\Phi_a(f) = 0 } . 
\label{}
\end{align}
Let us prove that $ D(\cG) = \tilde{D} $. 

Let $ f \in D(\cG) $. Let $ r>0 $ be fixed and set $ g=(r-\cG)f $. 
Then we have $ f=R_rg $. 
By \eqref{eq: Dynkin--}, \eqref{eq: rRrg-g--}, \eqref{eq: LR1--} 
and Proposition \ref{prop: L0R}, 
we have $ f,\cL f \in C_b([a,b)) $ and $ \Phi_a(f)=0 $. 
Hence we obtain $ D(\cG) \subset \tilde{D} $. 
By \eqref{eq: rRrg-g--}, we have 
\begin{align}
\cL f 
= \cL (R_rg) 
= rR_rg - g = \cG f . 
\label{}
\end{align}

Conversely, let $ f \in \tilde{D} $. 
Let $ r>0 $ be fixed and set $ g=(r-\cL)f \in C_b([a,b)) $. 
Then, by Proposition \ref{prop: RL}, we have 
\begin{align}
R^0_rg 
=& rR^0_rf - R^0_r \cL f 
\label{} \\
=& rR^0_rf - \cbra{ rR^0_rf - f + f(a) \frac{v_r}{v_r(a)} } 
\label{} \\
=& f - f(a) \frac{v_r}{v_r(a)} . 
\label{eq: Rr0g--}
\end{align}
Set $ h = R_rg-f $. 
By \eqref{eq: LR1--}, we have $ \Phi_a(R_rg)=0 $, 
and hence we have $ \Phi_a(h)=0 $. 
From \eqref{eq: Dynkin--} and \eqref{eq: Rr0g--}, it follows that 
\begin{align}
h = \cbra{ R_rg(a)-f(a) } \frac{v_r}{v_r(a)} . 
\label{}
\end{align}
By Lemma \ref{lem: Phiavr}, we have 
\begin{align}
0 = \Phi_a(h) 
= \cbra{ R_rg(a)-f(a) } \frac{\Phi_a(v_r)}{v_r(a)} 
= \cbra{ R_rg(a)-f(a) } \psi^{(b)}_a(r) . 
\label{}
\end{align}
By the condition \eqref{eq: pcond2}, we have $ \psi^{(b)}_a(r)>0 $, 
so that we obtain $ R_rg(a)-f(a)=0 $. 
This shows that $ h=0 $, which implies that $ f=R_rg $. 
Now we conclude that $ D(\cG) \supset \tilde{D} $, 
and thus the proof is complete. 
}

\subsection{The cases of {\boldmath $ 2^{\circ}) $} and {\boldmath $ 3^{\circ}). $}}

We prove Theorem \ref{thm: main-nonacc}. 

\Proof[Proof of Theorem \ref{thm: main-nonacc}]{
Suppose that $ a $ is accessible and $ b $ is not. 
Let $ g \in \cB_b = \cB_b([a,b]) $. 
In this case, we have the Dynkin formula 
\begin{align}
R_rg = R^0_rg + R_rg(a) \frac{v_r}{v_r(a)} 
\quad \text{on $ [a,b) $}. 
\label{eq: Dynkin-+}
\end{align}
Hence we have $ R_rg \in D(\cL) $ and we have the formula 
\begin{align}
\cL R_rg 
= r R_rg - g 
\quad \text{on $ [a,b) $}. 
\label{eq: rRrg-g-+}
\end{align}
In the same way as \eqref{eq: main--1}-\eqref{eq: main--2}, we have 
\begin{align}
\Phi_a(R_rg) 
=& n^{(b)}_a[\e^{-rT_b};T_b<\infty ] R_rg(b) 
\label{} \\
=& p_4(\{ b \}) R_rg(b) . 
\label{}
\end{align}
Thus we obtain the following: 
\begin{align}
\text{if $ g \in \cB_b $, we have $ R_rg \in D(\cL) $ and $ \Phi_a(R_rg)=p_4(\{ b \}) R_rg(b) $}. 
\label{eq: LR3}
\end{align}

(1) Suppose that $ b $ is entrance and irregular-for-itself. Set 
\begin{align}
\tilde{D} = 
\{ f \in D(\cL) : f,\cL f \in C_b([a,b]), \ \Phi_a(f)=p_4(\{ b \}) f(b) \} . 
\label{}
\end{align}
Let us prove that $ D(\cG) = \tilde{D} $. 

Let $ f \in D(\cG) $. Let $ r>0 $ be fixed and set $ g=(r-\cG)f $. 
Then we have $ f=R_rg $. 
By \eqref{eq: Dynkin-+}, \eqref{eq: rRrg-g-+}, \eqref{eq: LR3} 
and Proposition \ref{prop: L0R}, 
we have $ f \in C_b([a,b)) $ with finite left limit $ f(b-) $, 
$ \cL f \in C_b([a,b]) $ and $ \Phi_a(f)=p_4(\{ b \}) f(b) $. 
Since $ P^{(a)}_b = P^{\rm stop}_b $, 
we see that the Dynkin formula \eqref{eq: Dynkin-+} holds also for $ x=b $. 
This shows that $ f(b-)=f(b) $, hence we obtain $ D(\cG) \subset \tilde{D} $. 
By \eqref{eq: rRrg-g-+}, relation \eqref{eq: G1} is now obvious. 

We suppose that $ f \in \tilde{D} $. 
Let $ r>0 $ be fixed and set $ g=(r-\cL)f \in C_b([a,b]) $. 
Then, by Proposition \ref{prop: RL}, we have 
\begin{align}
R^0_rg 
= f - f(a) \frac{v_r}{v_r(a)} . 
\label{eq: Rr0g-}
\end{align}
Set $ h = R_rg-f $. 
By \eqref{eq: LR3}, we have $ \Phi_a(R_rg)=p_4(\{ b \}) R_rg(b) $, 
and hence we have $ \Phi_a(h)=p_4(\{ b \}) h(b) $. 
From \eqref{eq: Dynkin-+} and \eqref{eq: Rr0g-}, it follows that 
\begin{align}
h = \cbra{ R_rg(a)-f(a) } \frac{v_r}{v_r(a)} . 
\label{eq: h3}
\end{align}
Hence, by Lemma \ref{lem: Phiavr} and by \eqref{eq: stop hitting}, we have 
\begin{align}
0 
=& \Phi_a(h) - p_4(\{ b \}) h(b) 
\label{} \\
=& \cbra{ R_rg(a)-f(a) } 
\cbra{ \frac{\Phi_a(v_r)}{v_r(a)} - p_4(\{ b \}) \frac{v_r(b)}{v_r(a)} } 
\label{} \\
=& \cbra{ R_rg(a)-f(a) } 
\cbra{ \psi^{(b)}_a(r) - p_4(\{ b \}) P^{\rm stop}_b[\e^{-r T_a}] } . 
\label{eq: Phiap4b}
\end{align}
Since 
\begin{align}
\lim_{r \to \infty } \psi^{(b)}_a(r) = \infty 
\label{}
\end{align}
by the assumption \eqref{eq: pcond2} and since 
\begin{align}
\lim_{r \to \infty } P^{\rm stop}_b[\e^{-r T_a}] = 0 , 
\label{}
\end{align}
we see that $ \psi^{(b)}_a(r_0) - p_4(\{ b \}) P^{\rm stop}_b[\e^{-r_0 T_a}] > 0 $ 
for some $ r_0>0 $. 
Hence, by \eqref{eq: Phiap4b} we obtain $ R_{r_0}g(a)-f(a) = 0 $ 
and by \eqref{eq: h3} we obtain $ f=R_{r_0}g $. 
Now we conclude that $ D(\cG) \supset \tilde{D} $. 

(2) Suppose that $ b $ is natural or [entrance and regular-for-itself]. Set 
\begin{align}
\tilde{D} = 
\{ f \in D(\cL) : f,\cL f \in C_b([a,b]), \ \Phi_a(f)=p_4(\{ b \}) f(b) , \ \Phi_b(f)=0 \} . 
\label{}
\end{align}
Let us prove that $ D(\cG) = \tilde{D} $. 

Let $ g \in D(\cG) $. 
Let $ r>0 $ be fixed and set $ g=(r-\cG)f $. 
In the same way as (1), we can prove that $ f=R_rg $ and that 
$ f \in C_b([a,b)) $ with finite left limit $ f(b-) $, 
$ \cL f \in C_b([a,b]) $ and $ \Phi_a(f)=p_4(\{ b \}) f(b) $. 
For any $ r>0 $, we can find $ g \in C_b([a,b]) $ such that $ f = R_rg $. 
Using Theorem \ref{thm: main-1} where the roles of $ a $ and $ b $ are switched, 
and using the Dynkin formula \eqref{eq: Dynkin-+}, 
we have 
\begin{align}
\psi^{(a)}_b(r) f(b) 
=& q_3g(b) + \int_{(a,b)} q_4(\d x) R^0_rg(x) + q_4(\{ a \}) R_rg(a) 
\label{} \\
=& q_3g(b) + q_4[f] . 
\label{}
\end{align}
Noting that $ g = (r-\cG)f $, that $ g(b)=g(b-)=rf(b-)-\cL f(b) $, and that 
\begin{align}
\psi^{(a)}_b(r) = q_3 r + n^{(a)}_b[1-\e^{-r T_b}] = q_1 + q_3 r + q_4([a,b)) , 
\label{}
\end{align}
we have 
\begin{align}
\{ q_1 + q_4([a,b)) \} f(b) + q_3 r \{ f(b)-f(b-) \} + q_3 \cL f(b) = q_4[f] . 
\label{}
\end{align}
This shows that 
\begin{align}
\Phi_b(f) + q_3 r \{ f(b)-f(b-) \} = 0 . 
\label{}
\end{align}
Since $ q_3>0 $ and since $ r>0 $ is arbitrary, we obtain 
$ \Phi_b(f) = f(b)-f(b-) = 0 $. 
Hence we obtain $ D(\cG) \subset \tilde{D} $. 

Suppose that $ f \in \tilde{D} $. 
Let $ r>0 $ be fixed. 
Set $ g=(r-\cL)f \in C_b([a,b]) $. 
Then, in the same way as (1), we can prove that $ f = R_rg $ on $ [a,b] $, 
and hence we obtain $ D(\cG) \supset \tilde{D} $. 

The proof is now complete. 
}

The proof of Theorem \ref{thm: main-nonacc2} is quite similar 
to that of Theorem \ref{thm: main-nonacc}, and so we omit it.

\subsection{The case of {\boldmath $ 4^{\circ}). $}}

Now we prove Theorem \ref{thm: main}

\Proof[Proof of Theorem \ref{thm: main}]{
Suppose that both $ a $ and $ b $ are accessible. 
Let $ g \in \cB_b = \cB_b([a,b]) $. 
By the strong Markov property, we obtain the Dynkin formula: 
\begin{align}
R_rg = R^0_rg + R_rg(a) \frac{v_r}{v_r(a)} + R_rg(b) \frac{u_r}{u_r(b)} . 
\label{eq: Dynkin}
\end{align}
By Proposition \ref{prop: L0R}, we have $ R_rg \in D(\cL) $, and we have 
\begin{align}
\cL R_rg 
=& rR^0_rg-g + R_rg(a) \frac{rv_r}{v_r(a)} + R_rg(b) \frac{ru_r}{u_r(b)} 
\label{} \\
=& r R_rg - g . 
\label{eq: rRrg-g}
\end{align}
Using Lemmas \ref{lem: Phiavr} and \ref{lem: main4} 
and then using Theorem \ref{thm: main-1}, we have 
\begin{align}
\Phi_a(R_rg) 
=& \Phi_a(R^0_rg) + R_rg(a) \frac{\Phi_a(v_r)}{v_r(a)} + R_rg(b) \frac{\Phi_a(u_r)}{u_r(b)} 
\label{} \\
=& - N^{(b)}_{a,r}(g) + R_rg(a) \psi^{(b)}_a(r) 
+ R_rg(b) \cbra{ -n^{(b)}_a[\e^{-rT_b};T_b<\infty ] } 
\label{} \\
=& 0 . 
\label{}
\end{align}
Replacing the roles of $ a $ and $ b $, we obtain $ \Phi_b(R_rg) = 0 $. 
Thus we obtain the following: 
\begin{align}
\text{if $ g \in \cB_b $, we have $ R_rg \in D(\cL) $ and $ \Phi_a(R_rg)=\Phi_b(R_rg)=0 $}. 
\label{eq: LR1}
\end{align}

Let $ f \in D(\cG) $. Let $ r>0 $ be fixed and set $ g=(r-\cG)f $. 
Then we have $ f=R_rg $. 
By \eqref{eq: LR1}, we have 
$ \Phi_a(f)=\Phi_b(f)=0 $. 
Hence we see that $ D(\cG) $ 
is contained in the right-hand side of \eqref{eq: DG domain}. 
By \eqref{eq: rRrg-g}, we have 
\begin{align}
\cL f 
= \cL (R_rg) 
= rR_rg - g = \cG f . 
\label{}
\end{align}

Conversely, 
let $ f \in D(\cL) $ such that $ f,\cL f \in C_b([a,b]) $ 
and suppose that $ \Phi_a(f)=\Phi_b(f)=0 $. 
Let $ r>0 $ be fixed and set $ g=(r-\cL)f $. 
Then, by Proposition \ref{prop: RL}, we have 
\begin{align}
R^0_rg 
=& rR^0_rf - R^0_r \cL f 
\label{} \\
=& rR^0_rf - \cbra{ rR^0_rf - f + f(a) \frac{v_r}{v_r(a)} + f(b) \frac{u_r}{u_r(b)} } 
\label{} \\
=& f - f(a) \frac{v_r}{v_r(a)} - f(b) \frac{u_r}{u_r(b)} . 
\label{eq: Rr0g}
\end{align}
Set $ h = R_rg-f $. 
By \eqref{eq: LR1}, we have $ \Phi_a(R_rg)=\Phi_b(R_rg)=0 $, 
and hence we have $ \Phi_a(h)=\Phi_b(h)=0 $. 
From \eqref{eq: Dynkin} and \eqref{eq: Rr0g}, it follows that 
\begin{align}
h = \cbra{ R_rg(a)-f(a) } \frac{v_r}{v_r(a)} + \cbra{ R_rg(b)-f(b) } \frac{u_r}{u_r(b)} . 
\label{}
\end{align}
Since $ \Phi_a(h)=\Phi_b(h)=0 $, we obtain 
\begin{align}
\pmat{0 \\ 0} = A \pmat{R_rg(a)-f(a) \\ R_rg(b)-f(b) } , 
\label{}
\end{align}
where $ A $ is the matrix defined in \eqref{eq: A matrix}. 
Hence, by Proposition \ref{prop: main2}, we obtain 
\begin{align}
R_rg(a)-f(a) = R_rg(b)-f(b) = 0 . 
\label{}
\end{align}
This shows that $ h=0 $, which implies that $ f=R_rg $. 
Now we conclude that 
the right-hand side of \eqref{eq: DG domain} is contained in $ D(\cG) $, 
and thus the proof is complete. 
}

{\bf Acknowledgements.} 
The author expresses his sincere thanks to Professor Masatoshi Fukushima 
for his hearty encouragement and valuable comments. 
The author also thanks Professor Yukio Nagahata 
for his important question which motivated the author to obtain Theorem \ref{thm: main-nonacc}.

\end{document}